# *Une nouvelle démonstration d'irrationalité de racine carrée de 2 d'après les* Analytiques *d'Aristote*

**Résumé**. Pour rendre compte de la première démonstration d'existence d'une grandeur irrationnelle, les historiens des sciences et les commentateurs d'Aristote se réfèrent aux textes sur l'incommensurabilité de la diagonale qui se trouvent dans les *Premiers Analytiques*, en tant qu'ils sont les plus anciens sur la question. Les preuves usuelles proposées dérivent d'un même modèle qui se trouve à la fin du livre X des *Éléments* d'Euclide. Le problème est que ses conclusions, passant par la représentation des fractions comme rapport de deux entiers premiers entre eux i.e. la proposition VII.22 des *Éléments*, ne correspondent pas aux écrits aristotéliciens.

Dans cet article, nous proposons une nouvelle démonstration, conforme aux textes des *Analytiques*, fondés sur des résultats très anciens de la théorie du pair et de l'impair. Ne passant pas par la proposition VII.22, ni aucune autre propriété établie par l'absurde, cette irrationalité apparaît comme le premier résultat que l'on ne pouvait établir par une autre méthode.

L'importance de ce résultat, révélant un nouveau domaine mathématique, celui des grandeurs irrationnelles, rend compte de la centralité que cette forme de raisonnement acquiert alors, d'abord en mathématique, puis dans tout type de discours rationnel.

À partir des conséquences qui suivent de cette nouvelle démonstration, on peut interpréter très simplement la leçon sur les irrationnels du passage mathématique figurant dans le *Théétète* de Platon (147d-148b), ce que nous ferons dans un article à paraître dans un prochain volume.

***Abstract***. *To account for the first proof of existence of an irrational magnitude, historians of science as well as commentators of Aristotle refer to the texts on the incommensurability of the diagonal in* Prior Analytics*, since they are the most ancient on the subject.*
*The usual proofs suggested by the historians of science derive from a proposition found at the end of Book X of Euclid's* Elements*. But its conclusions, using the representation of fractions as a ratio of two integers relatively prime i.e. the proposition VII.22 of the* Elements*, do not match the Aristotelian texts.*
*In this article, we propose a new demonstration conformed to these texts. They are based on very old results of the odd/even theory. Since they use neither the proposition VII.22, nor any other result proved by a* reductio ad absurdum*, it seems to be the first result which was impossible to prove in another way.*
*The significance of this result, revealing a complete new territory in Mathematics, the field of irrational magnitudes, accounts for the centrality gained afterwards by this kind of reasoning, firstly in Mathematics, then in all forms of rational discourse.*
*From the consequences of this new proof, we can construe very simply the lecture on the irrationals in the mathematical text in Plato's* Theaetetus *(147d-148b). It will be done in an article to appear in a forthcoming volume.*

Salomon Ofman
Institut mathématique de Jussieu-Paris Rive Gauche
Université Paris 7
Histoire des sciences mathématiques
mail : salomon.ofman@imj-prg.cnrs.fr

# TABLE DES MATIÈRES

# Introduction.

La représentation géométrique des entiers par les Grecs anciens, permettait de les considérer comme une partie d'un ensemble plus vaste, celui des grandeurs linéaires. On obtenait ainsi la possibilité d'une représentation de grandeurs irrationnelles *i.e.* qui ne s'expriment pas comme un quotient de deux entiers. La question de l'origine de telles grandeurs irrationnelles ou « sans mesure commune » ou « incommensurables » (avec l'unité) est discutée. Néanmoins, à suivre la plus ancienne littérature qui nous est parvenue sur ce sujet, essentiellement les ouvrages de Platon et d'Aristote, il n'est guère douteux qu'elle suit de l'étude de la diagonale du carré, celle-ci étant « incommensurable » à son côté, ce que l'on note brièvement « incommensurabilité de la diagonale ». En langage moderne, cela s'exprime en disant que √2 est irrationnel[1].

Dans cet article, nous donnons une démonstration nouvelle de cette irrationalité, démonstration qui possède les propriétés suivantes.

(1) Elle utilise des résultats extrêmement anciens, connus sans doute des Mésopotamiens et certainement des Égyptiens de l'antiquité, il y a au moins 4000 ans, et probablement beaucoup plus[2].

L'intérêt d'une telle démonstration est sa neutralité chronologique. L'ancienneté des résultats nécessaires à la démonstration n'oppose aucune limite inférieure de datation, évitant ainsi les spéculations sur leur datation.

On pourrait donc s'interroger sur la tradition rapportant l'irrationalité aux Grecs anciens, et plus spécifiquement à Pythagore ou aux premiers pythagoriciens, c'est-à-dire autour du VI$^e$ siècle BCE[3]. C'est ce que nous ferons au dernier chapitre de la deuxième partie.

(2) Surtout, et c'est la raison principale qui a motivé notre recherche, elle est cohérente avec les textes (essentiellement d'Aristote, mais aussi de Platon). Pour reprendre une remarque de Maurice Caveing, elle donne sens à ce qui n'en a guère suivant les interprétations standard.

(3) Il n'est pas nécessaire de recourir, sous une forme ou sous une autre, à une proposition des *Éléments* d'Euclide sur la représentation des rationnels comme rapport de deux entiers premiers entre eux (la proposition 22 du livre VII). Non seulement cette proposition n'est pas évidente et suppose un développement important de la théorie des nombres, mais elle impose une limite inférieure à la datation de l'irrationalité de *√2*. Inversement, celle-ci force un point de vue chronologique sur cette proposition, et en conséquence sur l'ouvrage tout entier.

(4) On la retrouve comme en miroir au cœur de la dialectique platonicienne, en particulier dans sa stratégie définitionnelle, rendant compte de son importance chez le philosophe athénien.

(5)
Parce que, contrairement aux preuves usuelles, elle n'utilise aucun résultat démontré par l'impossible (ou l'absurde), l'incommensurabilité apparaît comme le premier résultat à n'avoir d'autre démonstration que par cette méthode de preuve. Ainsi, démonstration d'irrationalité et raisonnement par l'impossible entretiennent une relation très étroite, ce qui explique le rôle exemplaire des irrationnels, aussi bien chez Platon que chez Aristote.

---

[1] *Via* toutefois une certaine forme du théorème de Pythagore, cf. *infra* note 16. Rappelons que deux grandeurs *a* et *b* sont dites incommensurables ou sans commune mesure, si leur rapport *a/b* est irrationnel.

[2] Cela résulte du papyrus de Rhind, daté autour de 1600 BCE, copie d'un manuscrit plus ancien du Moyen empire vers 2 000 BCE. Toutefois, le texte semble être un manuel scolaire pour de jeunes garçons. Il ne saurait donner une vue d'ensemble des connaissances mathématiques auxquelles étaient parvenus, dans l'Égypte antique, ceux qui s'occupaient de ces questions, de même qu'un manuel de collège aujourd'hui ne saurait fournir un aperçu de l'ensemble des connaissances mathématiques contemporaines.

[3] Before Common Era.

(6) Ce n'est pas la particularisation d'un résultat plus large, contrairement aux démonstrations usuelles qui admettent des généralisations immédiates à tout entier non carré. Or, chez les auteurs anciens, l'irrationalité de √2 apparaît, sous la forme de l'incommensurabilité de la diagonale du carré, comme un résultat en soi. Et dans un texte du *Théétète* de Platon (147d-148b), portant sur l'irrationalité des racines d'entiers, le nombre *2* se trouve même explicitement distingué des autres entiers. À l'époque de Socrate, où se situe le récit, il n'y avait donc pas de démonstration unique, valide pour tous les nombres. Si celle visée par Aristote est bien la plus ancienne démonstration d'irrationalité, ce ne pouvait donc être une preuve générale, dont *2* serait un cas particulier.

Si même, comme le pensent certains commentateurs[4], il en allait différemment, et le récit du *Théétète* n'exposerait qu'une partie des résultats auxquels étaient parvenus les mathématiciens de cette époque, les textes que nous avons exigent une démonstration distincte pour *2*.

(7) Enfin, à l'inverse des interprétations habituelles qui se fondent sur la démonstration standard, notre démonstration permet d'expliquer de manière cohérente ce texte du *Théétète*, concernant l'irrationalité de certains entiers[5].

Cet article comprend trois parties.

Dans la première, qui est aussi la plus longue, on rappelle (chap. I-III) les diverses preuves classiques qui ont été proposées et qui dérivent d'un même modèle, que l'on trouve dans les *Éléments* d'Euclide. Nous leur avons donné le nom général de démonstration(s) standard. Le lecteur pressé, ou connaissant déjà les preuves classiques, pourra donc, au moins dans un premier temps, passer directement à la suite. Dans les deux derniers chapitres, nous confrontons tout d'abord ces démonstrations à la théorie de la connaissance aristotélicienne des *Analytiques*, puis aux textes aristotéliciens qu'elles sont censées expliquer.

La deuxième partie expose une nouvelle démonstration vérifiant les principales propriétés listées ci-dessus.

Dans la dernière partie, nous concluons sur les conséquences qui suivent de cette démonstration pour la méthode du raisonnement par l'impossible, et sa relation à la question de l'irrationalité. L'importance de ce sujet tient à la place qu'il occupe au carrefour entre histoire, mathématique et philosophie.

Pour les citations en langue étrangère, nous avons utilisé, autant que possible, les traductions existantes. Lorsqu'aucune n'était facilement disponible, nous avons traduit les textes, en donnant généralement l'original en note.

---

[4] Ainsi, quoique selon des schémas différents, Szabó 1977, p. 92-93, ou Caveing 1998, p. 162.

[5] Cela sera explicité dans un article à paraître dans un prochain numéro.

# Première Partie

## I. Présentation

L'irrationalité en mathématique intéresse non seulement les historiens des sciences, mais également les philosophes. Ainsi Aristote y fait-il plusieurs allusions dans les *Premiers analytiques*, sans en expliciter la démonstration. Par contre, il la commente de la manière suivante :

> *On prouve, par exemple, l'incommensurabilité de la diagonale, par cette raison que les nombres impairs deviendraient égaux aux nombres pairs, si on posait la diagonale commensurable ; on tire alors la conclusion que les nombres impairs deviennent égaux aux nombres pairs, et on prouve hypothétiquement l'incommensurabilité de la diagonale par ce qu'une conclusion fausse (ψεῦδος) découle de la proposition contradictoire. Car tel est, avons-nous dit, le raisonnement par l'absurde (διὰ τὸ ἀδυνάτου) : il consiste à prouver l'impossibilité (ἀδύνατον) d'une chose au moyen de l'hypothèse concédée au début. (An. pr. I, 23, 41a26-32, trad. Tricot.)*

Assertion qu'il reprend en la résumant un peu plus loin (I, 44, 50a37-38).

La démonstration dont il est question ici, est donc un raisonnement « par l'impossible » (ou « par l'absurde »)[6].

L'importance des textes aristotéliciens est considérable pour les historiens des mathématiques[7]. Ils permettent en particulier de donner une limite haute[8] à la connaissance des irrationnels par les mathématiciens grecs et plus particulièrement à la preuve visée par Aristote[9].

---

[6] Aristote emploie tantôt la formule « réduction à l'impossible » (*An. pr.* I, 7, 29b5), tantôt « preuve par l'impossible » (*ibid.* I, 21, 29b32), tantôt « preuve conduisant à l'impossible » (*An. post.* I, 24, 85a16).

[7] « Les travaux d'Aristote sont de la plus grande importance pour l'histoire des mathématiques et particulièrement des *Éléments*. Sa date (384-322/1) vient juste avant celle d'Euclide, en sorte que l'on peut tirer une assez bonne conclusion sur les innovations dues à Euclide lui-même, à partir des différences entre ses énoncés de choses correspondant à ce que nous trouvons dans Euclide, et les énoncés propres d'Euclide. » ('*The works of Aristotle are of the greatest importance to the history of mathematics and particularly of the Elements. His date (384-322/1) comes just before that of Euclid, so that from the differences between his statement of things corresponding to what we find in Euclid and Euclid's own we can draw a fair inference as to the innovations which were due to Euclid himself.*') (Heath 1981, p. 335).

[8] C'est-à-dire, la datation doit être antérieure à cette limite.

[9] Mais sans donner aucune limite basse. Ainsi que le remarquent M. Caveing ([CAV3], III, p. 132, note 6) et W. Knorr, il n'est aucune raison décisive pour la faire remonter au temps d'Aristote ou même de Platon. Notre analyse tendrait, au contraire, à montrer que le seul obstacle mathématique à surmonter, était une certaine forme de connaissance du théorème de Pythagore, même partielle et sans démonstration rigoureuse. Mais il en est un autre, qu'on pourrait qualifier de 'métamathématique', le plus difficile sans doute car concernant la forme même du raisonnement ( cf. *infra*, deuxième partie chap. VI, et troisième partie).

## II. Les démonstrations de l'irrationalité de 'racine carrée de *2*'

Par commodité pour le lecteur, nous utiliserons un langage moderne pour exposer les démonstrations.

Lorsque cela paraît utile à la lisibilité des énoncés, on désignera la multiplication par '×', sinon, comme il est d'usage, le signe de multiplication sera omis. Ainsi la multiplication des nombres *a* et *b* sera notée soit *ab*, soit *a*×*b*.

**1. Démonstration standard élaborée à partir des *Éléments*.**

Authentique ou pas, la démonstration qui se trouve au livre X des *Éléments* d'Euclide (appendice 27) de l'édition de Heiberg[10], est très généralement considérée comme celle visée par Aristote dans les *Premiers Analytiques*[11].

Si certains, ainsi Árpád Szabó, y voient un texte ancien préservé par Euclide, pour d'autres, il s'agit d'une opération tardive faisant partie de la tradition du commentaire d'Aristote[12].

Dans ce dernier cas, l'objectif de ses auteurs était d'ajouter aux *Éléments* une proposition donnant une référence mathématique sérieuse à l'allusion d'Aristote dans les *Premiers analytiques*. Le recours fréquent, sans preuves textuelles, à des « traces fossilisées », s'expliquant par la volonté des commentateurs d'harmoniser les récits des auteurs grecs anciens.

C'est néanmoins la démonstration qui, selon la plupart des historiens, est la preuve originale[13], et c'est à elle que renvoient encore, pour explication, les traducteurs du texte aristotélicien[14], aussi bien que les auteurs d'ouvrages généraux[15].

Il s'agit de démontrer que la racine carrée de *2* n'est pas rationnelle, ou en termes plus proches des géomètres grecs[16], que la diagonale du carré est incommensurable à son côté[17].

Pour cela, on note *z* cette grandeur (*i.e.* $z^2 = 2$), et

**Hypothèse** : on suppose qu'il existe deux nombres (entiers) *p* et *q* en sorte que l'on ait $z = p/q$. On peut en outre supposer *p* et *q* premiers entre eux[18]. On va en déduire une contradiction, ce qui prouvera que *z* n'est pas rationnel.

**Résultat 1** : *p* est pair.

Les égalités :

$z^2 = 2$ et $z = p/q$

donnent (en élevant les deux membres de la seconde égalité au carré) :

$2 = (p/q)^2 = p^2/q^2$

---

[10] Pour l'énoncé original, on se reportera par exemple à Vitrac 1998, p. 414-415.
[11] *Cf.* par exemple Heath 1998, p. 22-23.
[12] *Cf.* par exemple, Vitrac 1998, p. 24, 412-414.
[13] Ainsi Heath 1981, p. 90-91.
[14] Par exemple, Tricot 1947 [= AR12] ou Tredennick 1938 [= AR12'].
[15] Ainsi Russell 1970, p. 86 ou Kline 1982, p. 104-105.
[16] « L'ontologie des irrationnels se place toujours dans un contexte géométrique. Les Grecs n'ont jamais introduit pour ceux-ci de notation ou de termes arithmétiques standard ; notre "√2" est appelé le "côté du carré 2" ou étudié *via* le rapport du côté et de la diagonale du carré. » ('*The ontology of irrationals is always in the geometric context. The Greeks never introduced standard arithmetical notation or dictions for them; our '√2' was called the 'side of the square 2' or was studied via the ratio of the side and diameter of the square.*') (Knorr 1975, p. 9-10.) Le sens des termes traduits respectivement par irrationnel et incommensurable (ἄρρητος et ἄλογος) n'est pas toujours très clair, et ces termes sont en outre utilisés, chez de nombreux auteurs, comme des synonymes.
[17] En effet, si *b* est la diagonale d'un carré de côté *a*, d'après le théorème de Pythagore, on a : $b^2 = 2a^2$ i.e. $b^2/a^2 = 2$. L'incommensurabilité (de *b* par rapport à *a*) équivaut au rapport *b/a* non rationnel, autrement dit à l'irrationalité de la racine carrée de *2*.
[18] Autrement dit, ces deux nombres n'ont pas de diviseur commun différent de l'unité. Ainsi *4* et *15* sont relativement premiers ; mais *6* et *9* ne le sont pas, puisque *3* est un diviseur à la fois de *6* et de *9*.

d'où
$$2q^2 = p^2,$$
et $p^2$ donc $p$ est pair (car le carré d'un pair est pair, celui d'un impair est impair).

**Résultat 2** : $q$ est impair.

$p$ et $q$ étant premiers entre eux, ils ne peuvent être simultanément pairs, d'où : $q$ est impair.

**Conséquence :**

Suivant le résultat 1, $p$ est pair i.e. $p = 2r$ où $r$ est un entier. On a donc :
$$2q^2 = p^2 = (2r)^2 = 4r^2,\ \text{d'où}\ q^2 = 2r^2\ \text{est pair},$$
et (comme ci-dessus, le carré d'un impair étant impair), on obtient :

$q$ est pair.

Mais d'après le résultat 2, le nombre $q$ est impair ; $q$ est donc à la fois pair et impair, ce qui est impossible.

En conséquence, l'hypothèse de départ est fausse, et $z$ n'est pas commensurable avec $a$.

Cette preuve utilise implicitement la division en deux parties, sans élément commun, des entiers en pairs et impairs, et pas seulement le fait qu'ils sont disjoints. Cette propriété peut se lire dans l'équivalence non explicitée entre « pairs et impairs » et nombres, que l'on trouve communément dans les textes platoniciens[19].

**2. La démonstration d'Alexandre commentant Aristote**

Dans la recherche des origines de l'irrationalité, il importe, souligne Knorr, d'utiliser des textes qui ne soient pas tirés exclusivement des *Éléments* ; ce à quoi, ajoute-t-il, n'ont pas songé les historiens précédents[20].

On trouve une telle démonstration dans le *Commentaire* des *Premiers Analytiques* par Alexandre d'Aphrodise (III[e] siècle)[21]. Elle comporte diverses maladresses relevées par les commentateurs[22]. Pour la commodité du lecteur, nous la transposons en langage moderne en conservant les symboles de la démonstration précédente.

Alexandre réfère tout d'abord au « théorème 4 du livre X » (qui, dans la numérotation moderne, est le troisième théorème du livre X des *Éléments*), pour formuler l'hypothèse de la démonstration par l'impossible, à savoir :

**Hypothèse :**

$z^2 = 2$ et $z$ est un rapport de deux nombres entiers.

Puis il choisit $p$ et $q$ des nombres les plus petits possibles en sorte que l'on ait $z = p/q$.

Il renvoie alors au livre VII des *Éléments* (plus précisément, sans la nommer, à la proposition 22) pour obtenir :

$p$ et $q$ sont premiers entre eux        (**'**).

Suivant la proposition VII.27, il en déduit :

$p^2$ et $q^2$ sont premiers entre eux        (***).

**Résultat 1'** : $p^2$ est pair.

Suivant l'hypothèse, on a :
$$z^2 = 2 = p^2/q^2,$$
d'où :
$$2q^2 = p^2,$$
et donc :

$p^2$ est pair.

**Résultat 1'bis** : $q^2$ est pair.

---

[19] Entre autres, *Charm.* 166a, *Euthyd.* 290c, *Gorg.* 451b, *Hipp. maj.* 367a, *Ion*, 537e, *Phaed.* 104a-b, *Phaedr.* 274c, *Pol.* 259e, *Prot.* 357a, *Resp.* 510c, 522e, *Theaet.* 198a.

[20] Knorr 1975, p. 228.

[21] Pour le texte original, on peut se reporter à Wallies 1883 [= CAG II 1], p. 260-261. On en trouvera une traduction dans Vitrac 1998, p. 412-413, et ce que Knorr appelle une « paraphrase » (en anglais) dans Knorr 1975, p. 228.

[22] Par exemple Knorr 1975, p. 228-229.

Sans référence à Euclide[23], Alexandre écrit : « parce que la moitié des nombres carrés divisibles en deux parties égales est aussi paire », on obtient :

$q^2$ est pair.

**Résultat 2'** : $p^2$ est impair

D'après (**'), $p$ et $q$ sont premiers entre eux. Or pour que deux nombres soient premiers entre eux « il faut que tous deux soient impairs, ou que l'un seul d'entre eux le soit ». D'après le résultat 1'bis, $q^2$ est pair, donc :

$p^2$ est impair.

**Conséquence des résultats 1' et 2'** :

*il a été démontré grâce à l'hypothèse qu'ils* [$p^2$ *et* $q^2$] *sont tous deux pairs ; les impairs sont donc égaux aux pairs*[24] *en supposant la diagonale commensurable avec le côté, ce qui est impossible.*

Aussi l'hypothèse de départ est fausse, *i.e.* il n'existe pas d'entiers $p$ et $q$ tels que l'on ait $z = p/q$ autrement dit $z$ n'est pas rationnel. La démonstration est ainsi achevée.

Contrairement à la précédente, la contradiction (être pair et impair), et plus généralement la démonstration, porte non pas sur $p$ (ou $q$) mais sur leurs carrés.

**3. Une autre démonstration dans les *Éléments*.**

Une autre démonstration, plus brève, se trouve dans certains manuscrits de l'ouvrage, ainsi que dans la version gréco-latine[25].

En utilisant les mêmes notations et la même hypothèse, on remarque directement que, suivant la proposition VII.22, $p$ et $q$ sont premiers entre eux. Puis

**Résultat 1''** : $q$ n'est pas l'unité.

En effet, si $q$ était l'unité, son carré le serait également, et on aurait alors $p^2 = 2$, ce qui est impossible, car *2* n'est pas un carré (d'un entier).

**Résultat 2''** : $p$ et $q$ ne sont pas premiers entre eux.

En effet, $z^2 = p^2/q^2 = 2$ i.e. $p^2 = 2q^2$, d'où (prop. VIII.14) $q$ est un diviseur de $p$, et d'après le résultat précédent $q$ est différent de l'unité. Les entiers $p$ et $q$ ne sont donc pas premiers entre eux.

**Conséquence** : cela est impossible (puisqu'on a choisi $p$ et $q$ premiers entre eux), et donc l'hypothèse, la rationalité de $z$, *i.e.* de la racine carrée de *2*, est fausse.

---

[23] Pour Knorr 1975, p. 229, l'absence de renvois aux *Éléments* dans la seconde partie de cette démonstration est un argument supplémentaire pour penser qu'elle a été conçue par Alexandre lui-même, et en tout cas que son origine n'est pas euclidienne.

[24] Pluriel, ici, très surprenant. Dans cette démonstration, de même que dans les précédentes, on infère, de l'hypothèse posée par l'« absurde », qu'il existe précisément un nombre, différent suivant les preuves, qui est à la fois pair et impair. Alexandre reprend toutefois ici les termes effectivement employés par Aristote dans les *Analytiques*.

[25] On en trouvera une traduction par exemple dans Vitrac 1998, p. 416.

## III. Sur l'authenticité des démonstrations figurant dans les *Éléments*.

Toutes ces démonstrations sont ce que nous appelons des démonstrations par l'absurde, et procèdent de la manière suivante :

pour établir une propriété (soit *H*), on la suppose fausse, c'est-à-dire que sa négation (*non H*) est vraie.

On montre alors que de cette propriété (*non H*) s'ensuit une conséquence fausse.

De là, on conclut à la fausseté de *non H* (supposée hypothétiquement vraie au début).

La négation de la propriété à démontrer (*non H*) étant fausse, *H* est vraie.

En termes plus simples, prouver une proposition revient à montrer que sa négation est fausse car elle entraîne une contradiction[26]. Pour montrer qu'une grandeur (la racine carrée de *2*) est irrationnelle, on suppose qu'elle ne l'est pas, *i.e.* qu'elle est rationnelle, et on en déduit une impossibilité.

Sur la démonstration qui se trouve à la fin des *Éléments* dans l'édition de Heiberg (cf. *supra* II.1), les historiens adoptent des points de vue différents. La plupart considèrent qu'elle n'est pas d'Euclide. Oskar Becker remarque ainsi que « les éditeurs modernes (et déjà E.F. August en 1829) ont été embarrassés par cette proposition et l'ont écartée du texte »[27].

Selon Becker, rejoint par Á. Szabó, il s'agit d'une très ancienne proposition bien antérieure à l'œuvre d'Euclide. C'est lui (ou ses éditeurs) qui l'aurai(en)t ajoutée à l'ouvrage, parce que jugée « digne d'être conservée » pour des raisons historiques et/ou l'importance qu'il(s) lui accordai(en)t en tant que témoignage essentiel des origines des mathématiques[28].

Cette preuve, considérée comme participant d'une théorie pré-euclidienne et rajoutée dans les *Éléments*, serait alors celle à laquelle ferait référence Aristote[29]. Selon certains partisans de cette thèse, qu'elle se trouve chez ce dernier témoigne de son ancienneté ; selon d'autres, c'est précisément l'évidence de son ancienneté qui montre qu'elle est la preuve visée par le philosophe. Comme très souvent dans les questions de datation, le danger de circularité est présent[30].

Ceux qui contestent cette argumentation pensent qu'il s'agit d'un ajout tardif, dans la tradition des commentateurs aristotéliciens[31], tel Alexandre, dont nous avons rapporté, au paragraphe précédent, l'explication du passage des *Premiers analytiques* sur l'incommensurabilité. Plusieurs raisons, notent-ils, rendent invraisemblable l'ancienneté de la démonstration ajoutée à la fin du livre X des *Éléments*. Ainsi, remarquent-ils, elle utilise d'autres propositions euclidiennes dont la datation ne saurait être très ancienne[32]. Ils pointent également sa similitude avec celle d'Alexandre[33], qui, plusieurs siècles après les *Éléments*, en reprend des propositions dans le cours de la preuve, sans la citer à aucun moment. La proximité des deux preuves et les références explicites à des propositions euclidiennes non citées dans celle d'Alexandre suggèrent que cette démonstration qu'on trouve dans les *Éléments* leur, en réalité, est postérieure. Et surtout, si elle existait au temps d'Alexandre, plusieurs siècles après les *Éléments*, celui-ci n'aurait pas eu besoin d'en redonner une démonstration et y aurait directement renvoyé[34]. Ils en infèrent qu'elle a été ajoutée très tardivement, et en tout cas après le III[e] siècle CE[35], donc bien après Euclide, et ne saurait être connue d'Aristote.

---

[26] Par un curieux retournement, selon la perspective de la logique moderne, au lieu d'une démonstration, cela devient généralement la définition même de la négation d'un énoncé. Cela montre les limites d'une approche trop strictement formaliste, c'est-à-dire suivant de trop près les principes de la (des ?) logique(s) contemporaine(s).

[27] « Die Lehre vom Geraden und Ungeraden… », *Quellen und Studien zur Geschichte der Mathematik*, 3 **(1935)**, 533-553, note 1 ; *cf.* aussi Szabó 1977, p. 233

[28] *Ibid*. p. 234.

[29] « La preuve traditionnelle du caractère irrationnel de certains nombres est indiquée par Aristote comme un exemple de *reductio ad absurdum* (…) C'est évidemment la preuve interpolée au dixième livre des *Éléments* d'Euclide » ('*The traditional proof of the irrational character of certain numbers is indicated by Aristotle as an example of* reductio ad absurdum *(…) This is evidently the proof interpolated in the tenth book of Euclid's Elements*') (Edward A. Maziarz & Thomas Greenwood, *Greek Mathematical Philosophy*, New York, 1968, p. 50).

[30] D'où la prudence qu'il est nécessaire d'observer concernant les tentatives de datation. Sur ce point, *cf.* Pellegrin 1993 [= AR15], p. 65-66, Dalimier-Pellegrin 2004 [= AR14], ou encore Aubenque 1962, p. 10-11.

[31] « Cela indique que ce théorème résultait d'une interpolation tardive, mais pré-théonienne, stimulée par le travail de commentaire sur Aristote. » (Knorr 1975, p. 229.) Knorr ajoute que la fin de démonstration d'Alexandre, dont il soupçonne celui-ci d'être l'auteur, est un « appendice » (*appendage*) d'une preuve qui n'utilise pas le pair et l'impair (*ibid*.).

Enfin, ces deux démonstrations se fondent sur la possibilité de représenter toute fraction par un rapport de deux entiers premiers entre eux. Or, seul est nécessaire qu'ils soient de parités différentes. On peut en conclure, contrairement à Szabó et à la plupart des commentateurs, qu'il est peu probable que la démonstration visée par Aristote soit l'une ou l'autre de ces preuves[36].

---

[32] *Cf.* par exemple Knorr 1975, p. 24-26.
[33] Cf. *supra* II.2.
[34] Alexandre « possédait une copie des *Éléments* qui était essentiellement de la même forme que celle des manuscrits dont nous disposons, et sa copie ne contenait pas de preuve de l'incommensurabilité du côté au diamètre basée sur les propriétés des nombres impairs et pairs ; il a tiré cette preuve de quelque autre source. » (Knorr 1975, p. 229).
[35] Common Era.
[36] Par exemple Vitrac 1998, p. 417

# IV. Démonstrations des *Éléments* et causalité

Le dernier argument est d'autant plus sérieux que, dans tous les cas, la pierre de touche permettant de décider de la plus ancienne démonstration d'irrationalité sera les textes des *Premiers analytiques* (auxquels on peut ajouter celui, très bref, de *Métaphysique*, 983a19 *sq.*). Or selon Aristote, l'exactitude d'une démonstration ne dépend pas tant de la rigueur formelle du raisonnement (*i.e.* l'exactitude logique des enchaînements), que de la relation entre prémisse et conclusion. La prémisse doit avoir la généralité (« l'universalité ») qui convient exactement à la conclusion (*An. Post.* I, 2, 71b9-13).

Il est vrai que le rapprochement entre les démonstrations précédentes et celle dont parle Aristote paraît naturel.

D'une part, Aristote, insiste sur le type de preuve utilisé, celui par l'impossible. Parce que d'une hypothèse on déduit une impossibilité, il s'ensuit que la négation de cette hypothèse est vraie. Telle est bien la structure démonstrative euclidienne.

En outre, cette impossibilité est liée au problème d'une violation de la théorie du pair et de l'impair. De l'hypothèse posée découle en effet cette absurdité, qu'un même entier est à la fois impair et pair. Et la conclusion à laquelle aboutit la preuve rapportée par Aristote est bien une impossibilité relative au pair et à l'impair[37].

Toutefois, nous avons exposé au paragraphe précédent des arguments qui rendent un tel rapprochement douteux. Nous allons voir que la théorie aristotélicienne de la démonstration elle-même, telle qu'on la trouve dans les *Analytiques*, est incompatible avec cette présentation.

En effet, une caractéristique des démonstrations précédentes est qu'elles sont trop générales pour le résultat qu'il s'agit d'établir.

Or, selon Aristote, la connaissance ne consiste pas à énoncer une propriété vraie, ce n'est pas même son caractère principal. Il s'agit encore (*An. Post.* I, 24) d'en déterminer la *cause*, afin, dit-il, d'en obtenir une connaissance véritable ou universelle (καθόλου). Du point de vue démonstratif, cela se traduit par un raisonnement qui doit lui être précisément adapté. Il ne peut être ni trop général ni trop particulier. Dans l'un ou l'autre cas, affirme-t-il, la conclusion est incorrecte, même si, considérée « en elle-même », elle est vraie[38].

De même, la définition doit inclure la cause de son objet. Définir une chose consiste à la connaître, c'est-à-dire à en connaître l'essence[39].

Dans un autre ouvrage, il souligne encore :

*Il ne faut pas en effet se contenter, dans l'énoncé d'une définition, d'exprimer un fait, comme c'est le cas dans la plupart des définitions : il faut aussi que la cause* (αἰτίαν) *y soit présente et rendue manifeste* (ἐνυπάρχειν καὶ ἐμφαίνεσθαι)[40].

L'exemple qui suit ce texte, une quadrature, est, là encore, lié au rapport entre hypoténuse et côté d'un triangle rectangle. Dans le cadre de la géométrie grecque ancienne, il s'agit d'obtenir un carré (de surface) égale (à celle) d'un triangle ou d'un rectangle[41].

Pourtant, selon Aristote, cela n'est pas la réponse scientifique au problème. En effet, continue-t-il, de cette manière on obtient seulement la conclusion et non la définition cherchée (de la quadrature). Celle-ci est la « découverte d'une moyenne proportionnelle »[42], car c'est alors seulement, que « la cause de la chose » est exprimée ((τοῦ πράγματος λέγει τὸ αἴτιον)[43].

---

[37] Cf. *supra* I.
[38] En termes modernes, si dans le contexte considéré, sa valeur de vérité est vraie. On peut rapprocher cette théorie de la critique par Socrate de la science comme opinion vraie, à la fin du *Théétète*.
[39] Arist. *Metaph.* B, 2, 996b19-20.
[40] *De an.* II, 2, 413a13-15.
[41] C'est ainsi que procèdent Théétète et Socrate le Jeune dans le *Théétète* (148a-b) pour définir les racines carrées irrationnelles.

Dans ce cadre, la preuve de l'irrationalité de la racine carrée de *2* selon une méthode permettant d'établir aussi bien l'irrationalité d'autres nombres, est une connaissance « à la manière des sophistes »[44]. Elle ne rendrait, en effet, pas compte de la seule nature du nombre *2*.

Or les démonstrations précédentes ne vérifient pas cette condition aristotélicienne. Leur point commun est l'appel à des propositions des *Éléments* d'Euclide beaucoup plus générales que la conclusion obtenue, en contradiction avec la doctrine aristotélicienne de la science[45].

La dernière démonstration (*supra* II.3), certes la plus élégante, est également celle qui est la plus éloignée de la conception démonstrative aristotélicienne. Elle passe en effet, outre par la proposition VII.22, par la proposition VIII.14. En termes modernes, celle-ci établit une équivalence entre la division d'un nombre par un autre et celle de leurs carrés (i.e. $q$ divise $p$ équivaut à $q^2$ divise $p^2$).

De ces deux propositions, il s'ensuit immédiatement que :

La racine carrée de tout entier qui n'est pas un carré (d'entier) est irrationnelle.

Il suffit pour cela de remplacer, dans la preuve, *2* par un entier quelconque non carré. En voici une brève démonstration.

Soit $z$ la racine carrée de $n$ i.e. $z^2 = n$. D'après l'hypothèse par l'absurde (rationalité de $z$), on a :

$z = p/q$,

où $p$ et $q$ sont des entiers, qu'on peut supposer, suivant la proposition VII.22, premiers entre eux. En élevant les deux membres de l'égalité précédente au carré, on obtient :

$z^2 q^2 = n\, q^2 = p^2$.

**Résultat 1'' ** : $q$ n'est pas l'unité.

---

[42] La question soulevée permet de mieux saisir la notion de connaissance scientifique chez Aristote. En termes modernes, il affirme que la solution (scientifique) pour aboutir à la moyenne géométrique de deux nombres (soient *a* et *b*), ne consiste pas à construire un carré de surface égale à celle d'un triangle rectangle. Pourtant celle-ci donne un carré de côté *c*, et un triangle rectangle de côtés *2a* et *b*, en sorte que l'on a :

$c^2$ (surface du carré) = *(1/2)(2a)b* = *ab* (*cf.* figure ci-dessous) :

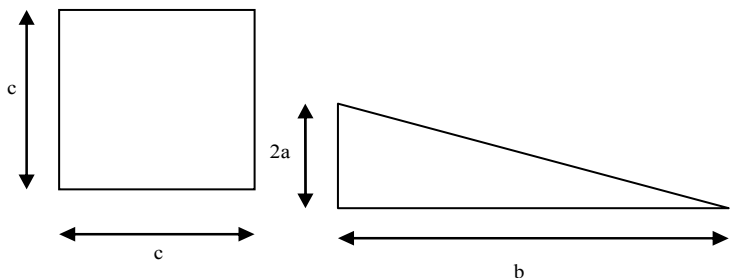

Mais selon le Stagirite, cela ne répond pas de manière satisfaisante au problème posé. On peut l'interpréter en comprenant que la construction géométrique fait apparaître l'existence d'une moyenne (le fait), sans en expliquer la cause (le pourquoi de cette existence).

[43] *De an.* II, 2, 413a15-20.

[44] *An. Post.* I, 2, 71b9-13.

[45] Bien que n'invoquant pas cette doctrine, Knorr remarquait déjà la distinction entre la démarche suivie dans ces démonstrations et le résultat recherché. Ainsi, visant la démonstration qu'on trouve à l'appendice du livre X des *Éléments* (cf. *supra* II.1), il écrit : « Un second témoignage de retouche dans la preuve est son utilisation de théorèmes trop généraux pour les conclusions voulues » ('*A second evidence of reworking in this proof is its use of theorems which are too general for its purposes.*') (Knorr 1975, p. 25) ; *ibid.* p. 229, concernant la démonstration d'Alexandre : « les principes invoqués (en particulier VII.27) sont beaucoup trop forts pour ce qui est demandé », ('*the principles invoked (in particular, VII, 27) are far stronger than required*') ; même critique portant cette fois contre la dernière démonstration (cf. *supra* II.3) qui donnerait l'irrationalité de la racine carrée de tout nombre qui n'est pas un carré (*op. cit.* p. 231).

Sinon (encore un raisonnement par l'impossible) son carré le serait également, et on aurait
$n q^2 = n = p^2$,

ce qui est impossible puisque (par hypothèse) $n$ n'est pas un carré (d'entier).

**Résultat 2''** : $p$ et $q$ ne sont pas premiers entre eux.

En effet, $p^2 = nq^2$ implique :

$q^2$ est un diviseur de $p^2$,

d'où $p^2$ et $q^2$ possèdent un facteur premier commun (à savoir $q^2$) qui, d'après le résultat 1'', est différent de l'unité, d'où (prop. VIII.14) :

$p$ et $q$ admettent un diviseur commun différent de l'unité. Ils ne sont donc pas premiers entre eux.

**Conséquence** : cela est impossible (puisqu'on a choisi $p$ et $q$ premiers entre eux). L'hypothèse est donc fausse, et la racine carrée de $n$ (à savoir $z$) n'est pas rationnelle.

On obtient ainsi non pas que « la racine carrée de *2* n'est pas rationnelle », mais la propriété beaucoup plus générale : la racine carrée de tout entier non carré (d'entier) n'est pas rationnelle.

Au sens d'Aristote, l'utiliser pour le seul cas particulier de la racine carrée de *2*, ne permet pas de remonter à la cause. Bien que le résultat soit vrai, c'est une vérité « sophistique » ou « accidentelle » (κατα συμβεβηκός)[46].

De toutes celles que nous avons considérées, la démonstration « usuelle » (*supra* II.1) se rapprocherait le plus d'une démonstration « minimale », c'est-à-dire utilisant les seules propriétés nécessaires à l'obtention de ce que l'on recherche.

Cependant, comme les autres, elle passe par la proposition VII.22, c'est-à-dire une propriété non de parité, mais de relative primalité[47]. Si la première concerne un entier particulier, l'autre porte sur un couple d'entiers[48]. Cette fois encore, le moyen n'est pas adéquat à la conclusion. En effet, pour l'obtenir, il est besoin non que toute fraction puisse être mise sous la forme d'un rapport de deux nombres premiers entre eux, mais seulement non simultanément pairs.

À considérer plus soigneusement les deux premières démonstrations, elles ne diffèrent d'ailleurs pas fondamentalement de la dernière. Car dans toutes les trois, la contradiction porte non pas tant sur le double attribut pair et impair d'un nombre indéterminé que sur le *couple* dont le rapport donne (sous l'hypothèse par l'absurde de sa rationalité) √2. En effet, les nombres qui le forment sont choisis relativement premiers, et on aboutit à ce qu'à la fois, ils sont relativement premiers et ne le sont pas. L'intervention de la parité étant en quelque sorte accidentelle.

Cela suit de l'appel à la proposition VII.22. Car la démonstration générale que nous avons donnée (pour l'irrationalité de la racine carrée de tout entier non carré d'entier) reste identique en utilisant, en place de la proposition VIII.14, la proposition VII.27[49]. Or celle-ci est une suite logique et une conséquence très simple de la première, ainsi que le montre la disposition de cette partie du livre VII, où elle apparaît comme une sorte de conclusion à un bloc ouvert par la proposition VII.22[50]. Il paraît donc improbable qu'un intervalle de temps important les ait séparées[51].

---

[46] L'utilisation de ce terme indique une opposition extrême à tout ce qui relève de la connaissance scientifique et à son caractère de nécessité. Le modèle d'une découverte accidentelle est celle où « creusant une fosse pour planter un arbre, on trouve un trésor » (*Metaph.* Δ, 30, 1025a15 ; *cf.* aussi E, 2, 1026b13-14).

[47] Cf. *supra* II.1, note 17.

[48] Cette proposition affirme que tout rapport est exprimable comme un rapport de deux entiers relativement premiers, ces entiers étant les plus petits nombres donnant ce rapport : ainsi le rapport 9/6 est égal au rapport 3/2.

[49] Celle-ci affirme que si deux entiers sont premiers entre eux, il en est de même de leurs carrés. De l'égalité $p^2 = nq^2$ où $q$ n'est pas l'unité (d'après le résultat 1''), on a $p^2$ et $q^2$ ne sont pas premiers entre eux, d'où (prop. VII.27), $p$ et $q$ ne le sont donc pas. On obtient donc le résultat 2'', et l'impossibilité puisque $p$ et $q$ ont été choisis premiers entre eux.

Selon ses partisans, un argument allant dans le sens de l'ancienneté de la démonstration « usuelle » serait sa référence à la proposition 23 du livre IX des *Éléments*[52]. Celle-ci suit de la division des nombres en pairs et impairs. Son ancienneté présumée, et plus généralement celle comprise entre les propositions 21 et 34 (voire 36, *cf.* Szabó 1977, p. 282) est un témoignage essentiel invoqué à l'appui de leur thèse par ceux qui pensent qu'elle est la preuve visée par Aristote[53].

Les autres références au texte d'Euclide seraient des modifications de la preuve initiale, une sorte de commodité offerte au lecteur qui avait facilement accès à cet ouvrage. Elles ne sauraient donc être considérées comme réfutant leur interprétation.

Mais ce raisonnement *ad hoc* peut être retourné, et le texte sera, au contraire, considéré comme un écrit tardif, où des éléments anciens ont été retravaillés pour les rendre conformes au modèle démonstratif euclidien, détruisant ainsi la structure propre du raisonnement primitif.

Knorr plaide pour une telle interprétation, se fondant sur la distance entre le moyen (théorèmes généraux euclidiens) et la fin (l'irrationalité de la racine carrée de *2*). Il ajoute que ce n'est qu'à partir de l'époque euclidienne qu'une séparation rigoureuse a lieu entre nombres et grandeurs, alors que l'auteur de cette démonstration distingue avec beaucoup de soin entre eux. L'approche ancienne aurait été plus « naïve »[54]. Il propose donc une reconstruction de la preuve pythagoricienne « purifiée de ces traits anachroniques », à partir de textes du *Ménon* de Platon.

Si l'on admet cette approche, la critique « interne » (*i.e.* provenant de la théorie aristotélicienne de la connaissance) que nous avons formulée à l'égard de la démonstration « usuelle » (*supra* II.1) n'est plus valable. Ainsi, de la représentation générale des fractions (d'entiers) par des nombres premiers entre eux, établie antérieurement, les mathématiciens grecs du temps d'Aristote auraient pu déduire, d'après un raisonnement sur le pair et l'impair, l'irrationalité de la seule racine carrée de *2*. Dès lors, n'eût-il pas été légitime de penser cette irrationalité comme conséquence de la nature particulière de ce nombre ? L'exemple serait sans doute infortuné, mais pas absolument contradictoire avec la théorie de la connaissance du Stagirite.

Toutefois, cette construction hypothétique ne s'accorde pas avec d'autres données textuelles. Car déjà au temps de Platon, et sans doute bien auparavant, il était connu que des racines carrées de nombres entiers différents de *2* étaient irrationnelles[55].

---

[50] H. Zeuthen remarque que le livre VII et le début du livre VIII forment un ensemble complet en lui-même (*cf.* Heath, p. 210).

[51] C'est pourquoi il est difficile de suivre Heath, lorsque, après avoir donné la démonstration standard d'irrationalité de racine carrée de *2*, il poursuit : « Cette preuve nous permet seulement de prouver l'incommensurabilité de la diagonale du carré avec son côté, ou de √*2* avec l'unité. Pour prouver l'incommensurabilité des côtés de carrés, dont l'un a *trois* fois la surface de l'autre, une procédure entièrement différente est nécessaire ; et on trouve en fait que, même un siècle après le temps de Pythagore, il était encore nécessaire d'utiliser des preuves *distinctes* (comme le passage du *Théétète* montre que Théodore le faisait) pour établir l'incommensurabilité avec l'unité de √*3*, √*5*, … jusqu'à √*17*. » ('*This proof only enables us to prove the incommensurability of the diagonal of a square with its side, or of √2 with unity. In order to prove the incommensurability of the sides of squares, one of which has* three *times the area of another, an entirely different procedure is necessary; and we find in fact that, even a century after Pythagoras's time, it was still necessary to use* separate *proofs (as the passage of the Theaetetus shows that Theodorus did) to establish the incommensurability with unity of √3, √5, … up to √17.*') ([Heath 1956, III, p. 2).

[52] Une somme d'un nombre impair d'impairs est encore impaire.

[53] Ainsi, Oskar Becker y voit le point de départ de constructions relevant d'une mathématique archaïque, ce qui, par ailleurs, soulève de sérieuses difficultés à propos de la cohérence de cet ensemble. *Cf.* par exemple Knorr 1975, p. 24 ; Vitrac 1998, p. 456-457.

[54] Knorr 1975, p. 25-26.

[55] *Cf.* par exemple Heath 1956, I, p. 413 et aussi III, p. 3 et 19.

Dans le *Théétète*, le personnage éponyme, un jeune athénien, rapporte à Socrate une leçon du mathématicien Théodore. Celui-ci montre que les segments formant les côtés des carrés de *3* et de *5* pieds (carrés) sont incommensurables avec le segment d'un pied, puis il continue avec d'autres nombres jusqu'à *17* (147d). Il est donc impossible qu'Aristote, ou les mathématiciens de son époque, aient pu penser que l'irrationalité était liée à la nature du nombre *2*, et de lui seul. Et de fait, l'exposé de Théodore exclut précisément ce nombre[56].

---

[56] Si le rapprochement avec la question considérée ici a souvent été fait, on n'en a guère tiré de conséquences pour l'interprétation du texte du *Théétète*. L'article sur ce passage du *Théétète*, à paraître dans un prochain numéro, rend compte de la leçon mathématique à partir de la démonstration que nous proposons ici (cf. *infra* deuxième partie) de l'irrationalité de la racine carrée de *2*.

## IV. Le texte d'Aristote : un pluriel bien singulier

Néanmoins la principale difficulté à laquelle se heurtent toutes les interprétations « standard » se trouve dans le témoignage aristotélicien lui-même. En effet, la conclusion (impossible) à laquelle on parvient dans le texte, est non pas qu'*un* nombre impair serait égal à *un* nombre pair ; moins encore qu'un nombre serait à la fois pair et impair, ce à quoi aboutissent, pourtant, toutes ces démonstrations.

L'assertion d'Aristote est beaucoup plus forte, puisqu'il conclut que, sous l'hypothèse par l'absurde :

*les nombres impairs deviendraient pairs*.

Autrement dit il n'y aurait plus de nombres impairs. Et en conséquence, mais seulement dans un deuxième temps, plus de nombres du tout. Aristote le répète deux fois en des termes identiques :

γίγνεσθαι τὰ περιττὰ (*les* impairs) ἴσα (égaux) τοῖς ἀρτίοις [*aux* pairs] : les <nombres> impairs deviennent égaux aux <nombres> pairs *An. pr.* I, 23, 41a27),

et une ligne plus bas :

ἴσα γίγνεσθαι τὰ περιττὰ τοῖς ἀρτίοις (41a28).

Un peu plus loin encore, il résume, utilisant toujours les mêmes termes, en se plaçant sur le mode du repos au lieu du mouvement (εἶναι au lieu de γίγνεσθαι), puisque, alors, la démonstration est achevée :

τὰ περιττὰ ἴσα εἶναι τοῖς ἀρτίοις (« les nombres impairs sont égaux aux nombres pairs » (*An. pr.* I, 44, 50a38).

Dans ces trois textes, la contradiction n'est jamais présentée sous la forme d'un double prédicat associé à un nombre, à la fois pair et impair, mais d'une relation d'égalité (*via* le terme ἴσα) entre la totalité des nombres impairs et celle des nombres pairs. Cela semble pourtant largement ignoré par les commentateurs. Ainsi, Becker-Hofmann écrivent dans leur *Histoire des mathématiques*[57] :

> *On trouve en outre dans Aristote (*Premiers ana.*, I, 23, 41a26-27) une remarque d'après laquelle si le côté du carré était commensurable avec la diagonale,* ce nombre *devrait être pair et impair.* (nous soulignons).

La position générale, avec diverses variantes, est bien représentée par celle de Thomas Heath. Il considère que les démonstrations d'Alexandre (démonstration II.2), celle que l'on trouve dans les *Éléments* (démonstration II.1), et celle à laquelle Aristote, par deux fois, fait allusion, sont « pratiquement de la même forme »[58],cette dernière n'étant rien d'autre que celle que l'on trouve encore dans les ouvrages modernes[59]. Celle proposée par l'historien anglais est ce que nous avons appelé la « preuve usuelle » (démonstration II.1). Elle conclut effectivement qu'un nombre (noté *q* dans cette démonstration) est à la fois impair et pair, « ce qui est impossible »[60].

---

[57] Becker & Hofmann 1956.

[58] '*It is given in practically the same form*' (Heath 1998, p. 22).

[59] « L'allusion concerne clairement la preuve bien connue de l'incommensurabilité de la diagonale que l'on trouve dans nos manuels scolaires. » ('*The allusion is clearly to the well-known proof of incommensurability of the diagonal which appears in our text-books.*') (*Loc. cit.*)

[60] *Ibid.*, p. 23. Pourtant, lorsque dans sa traduction commentée des *Éléments* d'Euclide, il traduit lui-même le texte des *Analytiques*, Heath écrit bien : « si on suppose [la diagonale (du carré) commensurable], il s'ensuivra que les (nombres) impairs seront égaux aux (nombres) pairs. » ('*if it be assumed [the diagonal (of the square) commensurable], it will follow that odd (numbers) are equal to even (numbers).*') (Heath 1956, I, IX.6, p. 136). Maurice Caveing est l'un des rares à remarquer qu'il y ait là une difficulté. Selon lui, elle résulterait toutefois d'une mauvaise traduction (Caveing 1998, III, note 7 p. 132). Nous pensons au contraire que celle-ci qui est celle retenue par tous les traducteurs des *Analytiques*, rend bien le sens

Dans un autre passage du *Théétète* que celui déjà mentionné, Socrate exprime la même idée que celle des *Analytiques*. Il s'agit d'exemples d'absurdités si démesurées que, la raison étant endormie, l'esprit se refuse encore à les envisager. Même en rêve, dit-il, il est impossible de jamais penser que :

τὰ περιττὰ ἄρτιά ἐστιν (« les impairs sont pairs », 190b).

L'énormité d'une telle supposition est soulignée encore par sa relation à deux autres impossibilités, « se persuader que nécessairement le bœuf est un cheval, ou que deux [est] un » (190c).

Aristote insiste, lui aussi, sur l'évidence de la conclusion. À l'inverse de la plupart des raisonnements, non mathématiques, qui demandent l'accord des interlocuteurs sur certains points pour aboutir à leur fin, rien de tel ici. Elle suit du raisonnement sans qu'il y ait place pour le moindre doute :

*La différence avec les syllogismes hypothétiques précédents, c'est que, dans les premiers, une convention préalable est nécessaire pour entraîner l'assentiment à la conclusion (...) par contre, dans le cas présent, même en l'absence d'une convention préalable, on donne son assentiment en raison de l'évidence de l'erreur*[61].

Car l'unité étant ce par quoi le nombre est défini, à savoir une « pluralité » d'unités, si deux était un, tous les nombres seraient pairs, et on retrouverait cette absurdité des impairs devenant pairs, évidence telle qu'aucune argumentation n'est nécessaire.

En outre, la structure même de la phrase (τὰ περιττὰ ἄρτιά ἐστιν, « les impairs sont pairs ») est importante. Il s'agit d'attribuer aux « impairs » (τὰ περιττά) un attribut, celui d'être pairs. La conclusion s'oppose à une simple conjonction, que ce soit dans « un pair égal à un impair » ou « un nombre à la fois pair et impair ». Celle-ci est en effet symétrique. Car c'est la même chose qu'un nombre soit « pair et impair » ou « impair et pair » ; et encore, que « un pair est égal à un impair » ou « un impair est égal à un pair ». Les deux formulations se retrouvent d'ailleurs chez les commentateurs rapportant le texte aristotélicien.

Au contraire, la phrase considérée est asymétrique. Elle n'indique pas, du moins directement, une identification entre les pairs et les impairs, mais plutôt que les impairs seraient des pairs.

Il est vrai que dans ce texte, c'est le terme d'égalité (ἴσα) qui est utilisé par Aristote.

Un premier sens peut être immédiatement exclu, celui d'égalité numérique (ainsi dans l'expression ἴσας ναῦς, « des navires en nombre égal »[62]). Outre que la construction est différente, les ensembles considérés ici sont infinis. Or un argument d'Aristote pour rejeter le caractère numérique de l'infini, est précisément que tout nombre doit être ou pair ou impair, l'infini n'étant ni l'un ni l'autre[63].

Une autre possibilité consisterait à identifier les deux multitudes. L'asymétrie des trois énoncés successifs, où la propriété porte toujours sur les nombres impairs, va toutefois à l'encontre d'une telle interprétation. Certes, l'égalité mathématique est une relation symétrique, si *a* = *b* alors *b* = *a* et inversement. Mais si Aristote reprend un résultat mathématique, il le fait dans un cadre philosophique[64]. Or l'égalité, en langage usuel, n'est pas une relation symétrique.

---

du texte aristotélicien, l'incohérence avec la démonstration retenue, à savoir celle tirée des *Éléments*, provenant de ce qu'Aristote en visait une autre.

[61] *An. pr.* I, 44, 50a32-39 (trad. Tricot).
[62] Xénophon, *Hell.* I, 6, 30 : ἐπετέτακτο δὲ Πρωτομάχῳ μὲν Λυσίας, ἔχων τὰς ἴσας ναῦς (« en soutien de Protomachos, Lysias, avec le même nombre de vaisseaux »). Xénophon décrit la position des flottes adverses lors de la bataille des Arginuses (406 BCE) entre Sparte et Athènes. La phrase commence par : « L'aile droite était tenue par Protomachos avec quinze vaisseaux ; à côté de lui, Thrasyllos avec quinze autres ; en soutien … »
[63] *Metaph.* M, 8, 1084a1-2.
[64] Qu'il importe pour lui de distinguer du cadre mathématique, *cf. Metaph.* A, 9, 992a35-b1.

Ainsi dans le *Gorgias* de Platon, Socrate comparant les politiques anciens et modernes, dit, en employant ce même terme ἴσος : οὗτοι δὲ ἀνεφάνησαν ἐξ ἴσου τοῖς νῦν ὄντες (517a2), traduit par A. Croiset (1923) de la manière suivante : « ceux-ci nous sont apparus comme valant les modernes » ; par M. Canto (1987) : « ces hommes étaient pareils à nos contemporains » ; ou encore par W. Lamb (1925) : « nous les avons trouvés être au même niveau que les nôtres au jour d'aujourd'hui » ('*these we have found to be on a par with ours of the present day*')[65]. Dans tous les cas, la comparaison porte sur le premier terme relativement au second. À gauche ce que l'on recherche, à droite l'étalon de comparaison[66]. Il en est de même en français usuel, ainsi lorsqu'on parle d'égalité des droits[67].

On peut donc comprendre qu'Aristote n'utilise pas tant ἴσα au sens d'une relation d'égalité mathématique strictement identificatoire que d'une attribution. Le sens est bien que tous les impairs *deviendraient* pairs. D'où la dissymétrie systématique des trois énoncés aristotéliciens, où, sous l'hypothèse par l'impossible, ce sont toujours les *impairs* qui reçoivent l'attribution. L'ambiguïté de la phrase, tient à ce que, une fois *devenus* pairs, on a aussitôt une identification entre pairs et impairs[68]. Suivant cette analyse, Aristote et Platon réfèrent bien au même énoncé mathématique.

Quel que soit le sens précis que l'on prête ici à ἴσα, on aboutit à la disparition des impairs au profit des pairs. Situation qui est bien de la dernière absurdité, équivalant aux autres exemples donnés par Socrate dans le *Théétète*, impossible à imaginer, même en rêve.

Reste à considérer la démonstration visée par Aristote, lui permettant de conclure à l'inclusion des impairs dans les pairs (ou à leur égalité). Nous allons tout d'abord examiner les difficultés que posent les démonstrations précédentes vis-à-vis de cette question.

Celles que l'on trouve dans les divers manuscrits des *Éléments*, et que nous avons rapportées précédemment (*supra*, démonstrations II.1 à II.3) comportent quatre étapes que l'on rappelle brièvement.

**Hypothèse 1** : $z = p/q$.
**Propriété 1** : on choisit les entiers $p$ et $q$ premiers entre eux.
**Conséquence 1** : $2q^2 = p^2 = 4r^2$, d'où $q^2 = 2r^2$.
**Conséquence 2** : $p$ et $q$ sont pairs, ce qui contredit la propriété 1.

La propriété 1 est souvent posée comme une évidence donnée en passant, voire en note de bas de page. Si cela est cohérent avec le texte d'Aristote, lequel ne fait aucune allusion à un tel résultat, cela ne l'est guère dans le cadre d'une reconstitution de preuve mathématique. Cette propriété est en effet le point central de la démonstration.

La représentation de tout rapport par des nombres premiers entre eux présuppose déjà à la fois une arithmétique suffisamment élaborée et un raisonnement par l'impossible devenu une méthode de preuve standard en mathématique. Toutes les démonstrations précédentes réfèrent, en effet, implicitement ou explicitement, à la proposition VII.22 des *Éléments* :

---

[65] Voir respectivement Croiset 1923 [= PLA5'], Canto 1987 [= PLA5], Lamb 1925 [= PLA5"].

[66] Même asymétrie dans le schème de la division dialectique par dichotomie décrit par Socrate en *Phèdre*, 266a. Cette question de symétrie et d'asymétrie se trouve également au cœur du *Parménide*. Car si, comme l'affirme le personnage éponyme de Platon, la ressemblance est une relation nécessairement symétrique, c'est toute la théorie philosophique présentée par Socrate, et aussi bien celle de Platon, qui risque de s'effondrer (sur cette question extrêmement controversée, *cf.* p. ex. Goldschmidt 2003 [= GOL], p. 47-51 et en particulier p. 50 note 1). Mais ainsi que le remarque ce dernier, « personne, sans doute, n'aura l'idée de dire que le modèle ressemble à la copie ni qu'un père ressemble à son fils », si ce n'est, toutefois, « par boutade » (*op. cit.* p. 50).

[67] Il serait par exemple étrange que les amis des animaux exigent que « leurs droits soient égaux à ceux des hommes », afin de réduire ceux des hommes, et ramener ces derniers à l'état de bêtes. Il s'agit, *au contraire* en quelque sorte, de demander que les droits des premiers deviennent égaux à ceux des seconds.

[68] Dans la démonstration que nous proposons, sous l'hypothèse de la rationalité de la racine carrée de *2*, on conclut à cette impossibilité des impairs (en totalité) devenus pairs. L'identification s'ensuit alors, mais d'après une propriété supplémentaire, quoique évidente, à savoir que tout pair est somme d'un impair et de l'unité (cf. *infra* note 103). D'un point de vue négatif, cette inclusion impliquerait la parité de l'unité, contrairement à son indivisibilité définitionnelle.

*Les nombres les plus petits parmi ceux qui ont le même rapport qu'eux sont premiers entre eux*[69],
dont la preuve est obtenue par un raisonnement par l'impossible.

Si la preuve primitive d'incommensurabilité passait par cette proposition, le type de raisonnement par l'impossible ne serait donc pas apparu pour traiter de l'irrationalité, mais dans un cadre numérique (entier), pour démontrer ce résultat : tout rapport (d'entiers) est représentable par un rapport de deux entiers premiers entre eux[70].

L'exemple de l'irrationalité traité par Aristote perdrait ainsi de son intérêt historique et de sa pertinence, alors que, rappelle M. Caveing, suivant l'*Index Aristotelicus* d'Hermann Bonitz, il n'apparaît pas moins de vingt-six fois dans le *corpus* aristotélicien[71]. Le Stagirite ne l'aurait pas choisi en tant que résultat paradigmatique de ces démonstrations, mais arbitrairement. Il aurait même fait preuve d'une certaine incohérence, puisque l'exemple fondamental pour illustrer ce mode de raisonnement nécessiterait un résultat intermédiaire obtenu d'après ce même mode. Pour éviter toute circularité dans la présentation, c'est celui-ci qui aurait dû être choisi, le cas étant d'ailleurs plus simple. En procédant différemment, Aristote irait également à l'encontre de la théorie de la connaissance qu'il soutient dans les *Seconds analytiques*[72].

Quoi qu'il en soit, dans la reconstruction des démonstrations que nous avons considérées ci-dessus, l'un des deux raisonnements s'impose :
– l'incommensurabilité de la diagonale du carré à son côté doit être rapportée à la possibilité d'écrire un rapport numérique quelconque sous la forme de deux nombres premiers entre eux ;
– tout rapport numérique s'écrit comme un rapport de deux nombres premiers entre eux, sous peine de l'absurde commensurabilité de la diagonale au côté du carré.

Ces deux énoncés affirment la même chose, le premier sous une forme affirmative, le second sous une forme négative[73].

Pourtant, que ce soit chez Platon ou Aristote, on ne trouve aucune allusion à une construction imbriquant deux raisonnements par l'impossible, l'un établissant la propriété de représentation des rapports numériques par des entiers premiers entre eux, puis, de là, l'impossibilité de la rationalité de racine carrée de *2*.

Pour éviter ces difficultés, il va s'agir de rendre compte de l'irrationalité, suivant une démonstration qui n'utilise pas cette argumentation, et plus généralement ne passe pas par un résultat prouvé par l'impossible.

Enfin et surtout, il faut que la conclusion s'accorde avec celle des textes aristotéliciens, plutôt qu'avec ceux tirés des *Éléments* d'Euclide (authentiques ou pas) et de commentateurs tardifs.

---

[69] 'Οἱ ἐλάχιστοι ἀριθμοὶ τῶν τὸν αὐτὸν λόγον ἐχόντων αὐτοῖς πρῶτοι πρὸς ἀλλήλους εἰσίν', proposition passant par la définition VII.13 d'entiers premiers entre eux.
[70] Il faudrait alors revoir les datations concernant ce raisonnement en mathématique (*cf.* Szabó 1977, par exemple).
[71] Caveing 1998, t. 3, p. 133.
[72] Cf. *supra*, chapitre IV.
[73] À la manière dont Socrate, parlant du texte de Zénon, le met en relation avec la doctrine parménidienne : « D'une certaine façon, en effet, [Zénon] a écrit la *même chose* que [Parménide], mais en *retournant* l'argumentation (…) [Parménide] pose que l'univers est un (…) Zénon, lui, pose à l'inverse que les choses ne sont pas plusieurs (…) et chacun de votre côté (…) *vous dites la même chose*. » (Platon, *Parménide*, 128a-b, trad. Brisson [= PLA2] ; nous soulignons.)

# Deuxième Partie

## I. Dichotomie et décomposition pair/impair

Tout d'abord, la propriété 1 n'est pas nécessaire dans toute sa généralité : il suffit que $p$ et $q$ ne soient pas simultanément pairs[74]. Et il est vrai que l'appel à la proposition VII.22 des *Éléments* n'est pas indispensable[75]. La démonstration que nous allons donner repose en effet sur un raffinement de la très ancienne méthode du pair et de l'impair, qui, pour certains[76], précède même l'apparition de l'écriture. En tout état de cause, elle était bien connue des Égyptiens antiques qui utilisaient une représentation binaire (en parallèle à la représentation décimale) des nombres, particulièrement utile pour les multiplications ou les divisions[77].

Suivant les quelques papyrus qui nous sont parvenus de l'époque égyptienne antique, en particulier celui dit de Rhind (R37 et 38) daté entre 1500 et 1800 BCE, lui-même copie d'un texte plus ancien (autour de 2000 BCE), on peut considérer que la méthode du pair et de l'impair raffinée, que nous allons décrire ci-dessous date, au moins, de 4000 ans. Il s'agit de distinguer non pas seulement entre les nombres pairs (divisibles par *2*) et les impairs qui ne le sont pas, mais de continuer la division des pairs aussi loin que possible. Autrement dit, de considérer le nombre maximal de fois qu'un nombre peut être divisé par *2*, avant que l'on obtienne un nombre impair.

Cela s'écrit (en langage moderne) :
Pour tout entier $n$ pair[78], on a :
**$n = 2^h u$** où $h$ est le nombre de divisions possibles de $n$ par *2*, et $u$ est un impair.

On doit noter que dans cette écriture, $h$ (donc $u$) est bien défini dès que $n$ est donné. En effet, $h$ est **le** nombre maximal de divisions possible d'un entier par $2^{79}$. Ainsi :

$28 = 2 \times 14 = 2 \times 2 \times 7$ est divisible *2* fois par *2*, i.e. $h = 2$ ou encore $28 = 2^2 \times 7$,

$30 = 2 \times 15$ est divisible une seule fois pas *2*, i.e. $h = 1$ ou encore $30 = 2^1 \times 15$,

$32 = 2 \times 16 = 2 \times 2 \times 8 = 2 \times 2 \times 2 \times 4 = 2 \times 2 \times 2 \times 2 \times 2 = 2^5$, i.e. $h = 5$ ou encore $32 = 2^5$,

$40 = 2 \times 20 = 2 \times 2 \times 10 = 2 \times 2 \times 2 \times 5$, i.e. $h = 3$ ou encore $40 = 2^3 \times 5$.

L'écriture de tout entier pair sous la forme **$n = 2^h u$** (où $u$ est impair), bien que découlant naturellement de la distinction entre pair et impair, permet de découper l'ensemble des nombres de manière plus précise. En fait, il s'agit de répéter ce découpage, mais cette fois-ci sur les pairs.

Nous appellerons cette écriture la *décomposition pair/impair* de $n$. Et puisque dans la suite, l'entier $h$ joue un rôle essentiel, nous l'appellerons *degré de parité* de $n$.

---

[74] Cf. *supra* II.1. Bien qu'il reprenne la démonstration usuelle, Heath n'en souligne pas moins que (a) seul est nécessaire l'argument de parité/imparité, et non pas que $p$ et $q$ soient premiers entre eux ; (b) cette dernière propriété ($p$ et $q$ premiers entre eux) ne permet pas d'obtenir un résultat plus large (*i.e.* pour les racines carrées de nombres différents de *2*).

[75] C'est également ce que note le physicien-mathématicien contemporain Roger Penrose lorsqu'il propose une démonstration aussi simple que possible de l'irrationalité de racine carrée de *2* pour le lecteur contemporain (Penrose 2007, §3.1, p. 48 et note 2).

[76] Comme P.-H. Michel, cité dans Vuillemin 2001, p. 2.

[77] *Cf.* p. ex. Heath 1981, p. 29.

[78] Pour nous autres modernes, l'hypothèse $n$ pair est superflue, car, on considère alors que $h$ est égal à *0* et $2^0 = 1$. Dans le cadre grec ancien, on considérerait séparément les cas pair et impair.

[79] En termes modernes on parlerait de l'unicité de $h$ et de $u$ dans l'écriture de $n$ sous la forme $n = 2^h u$.

En notations modernes, si l'on pose :

**N** est l'ensemble des entiers, **I** l'ensemble des impairs et **P** l'ensemble des pairs,

**2I** l'ensemble des nombres impairs multipliés par *2* et **2P** l'ensemble des nombre pairs multipliés par *2*,

**4I** l'ensemble des impairs multipliés par *4* et **4P** celui des pairs multipliés par *4*,

et ainsi de suite,

on obtient le tableau de la figure 1 ci-dessous :

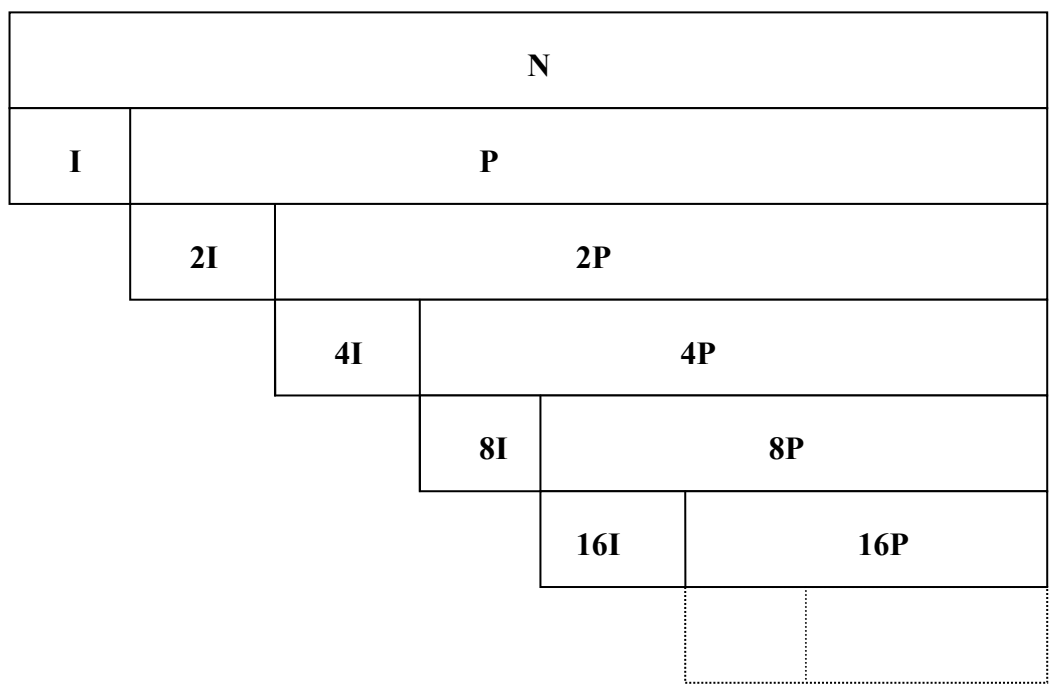

Fig. 1

Ce tableau, essentiellement dû à Nicomaque, mais que les historiens font remonter aux Pythagoriciens anciens, décrit précisément les pairs sous la forme d'une suite :

  *2(2n+1), 4(2n+1), 8(2n+1), 16(2n+1), …*[80].

C'est un cas très particulier d'écriture binaire. Dans le cadre opératoire égyptien, où elle était indispensable pour tout calcul, cette propriété allait de soi.

Dans l'esprit des mathématiciens de la Grèce antique, elle est équivalente à la finitude de la suite des divisions par *2* d'un entier quelconque. Cela signifie encore qu'il n'existe pas de descente infinie de divisions par *2* pour les entiers.

Comme dans le décryptage socratique de l'ouvrage de Zénon[81], on a un même raisonnement, présenté sous deux formes, l'une affirmative, l'autre négative. Un début de classification fondé sur cette propriété est élaboré au livre VII des *Éléments* d'Euclide[82]. Cette construction intervient également dans la recherche des nombres « parfaits » que l'on trouve à la proposition 36 du livre IX[83].

---

[80] Bertier 1978, I, X, 9, p. 67-68. Pour le lecteur intéressé, cela donne une démonstration très simple qu'il y a « autant » d'entiers que de fractions, c'est-à-dire qu'on peut numéroter la totalité des fractions par des entiers. En particulier, la totalité des nombres rationnels peut être numérotée par les seuls entiers (c'est ce genre de paradoxe qui a fait dire à Cantor, 'je le vois mais je ne le crois pas').

[81] Cf. *supra*, note 71.

[82] Cinq (ou six suivant que l'on tient ou non l'une d'elles pour apocryphe) des définitions par lesquelles débute ce livre concernent les divers types de nombres pairs, désignés comme étant pairement impairs ou pairement pairs. Cette construction conduit naturellement à la recherche de la parité maximale des entiers. Nombre de commentateurs antiques tardifs, dont Philopon, Nicomaque et Jamblique, donnent un sens différent à ces termes, afin d'obtenir une partition des entiers en fonction de leur parité, *via* la propriété 2. On va ainsi de ceux qui ont la forme d'une puissance de deux (« les plus pairs », pourrait-on dire) à leur extrême opposé, les impairs (*cf.* Heath 1956, II, p. 281-284).

Une telle suite des puissances de *2* se retrouve dans les plus anciens textes grecs non mathématiques, leur connaissance était donc certainement très antérieure. Déjà les raisonnements par dichotomie de Zénon se fondent sur elles et leurs inverses, *i.e.* les fractions de la forme *(1/2)$^n$*.

Et l'Athénien des *Lois* de Platon propose de fixer les châtiments pour récidives, de la manière suivante :

> *s'il est reconnu coupable d'avoir souillé par une occupation indigne le foyer qui est le sien et celui de ses ancêtres, il sera détaché de cette occupation en faisant un an de prison. Et s'il est repris, il fera deux ans de prison, et, à toute récidive il sera emprisonné de nouveau sans répit pour un temps chaque fois double du précédent*[84].

Plus généralement, la division platonicienne dans la recherche des définitions prend modèle sur cette suite[85].

Dans un texte critiquant le procédé dichotomique platonicien, Aristote remarque que celui-ci aboutit à ce que « les différences dernières seront au nombre de *4* ou ce sera quelque autre parmi les binaires successifs [en fait les puissances de *2*] »[86]. Il reproche ainsi à Platon d'affirmer que les seules divisions naturelles seraient des puissances de *2*, l'erreur de ce dernier étant d'oublier le facteur impair qui intervient lorsqu'on considère un nombre en général. Ce serait un autre exemple où, selon lui, un philosophe torture les faits pour imposer sa thèse[87].

C'est encore suivant une telle procédure que peut se comprendre l'énigmatique engendrement de tout nombre à partir de l'« Un » et de la « Dyade », attribué à Platon par Aristote, et que l'on retrouve fréquemment lorsque ce dernier conteste les thèses de l'Académie[88].

Un texte de la *Métaphysique* est particulièrement explicite. Après avoir noté, une fois encore, que la conception platonicienne du Nombre « procède de l'Un et de la Dyade indéfinie »[89], Aristote remarque qu'un nombre ne peut être infini, car dans ce cas, dit-il, il ne serait ni pair ni impair, or tout nombre est nécessairement l'un ou l'autre. Puis il ajoute :

---

[83] Un nombre est « parfait » (τέλειος), s'il est égal à la somme de tous ses diviseurs (sauf lui-même). Ainsi *6 = 1+2+3* est un nombre parfait. Ils permettent d'obtenir une décomposition de l'unité en fractions : ainsi, en divisant les deux membres de la dernière égalité par *6*, on obtient : *1 = 1/6 + 1/3 + 1/2*. La proposition d'Euclide montre que ce sujet est étroitement lié à ce qu'on appelle aujourd'hui les nombres *premiers de Mersenne* (*i.e.* les nombres premiers de la forme *($2^k$ - 1)* où *k* est un nombre premier). La recherche de ces nombres est difficile. Le dernier découvert, chronologiquement le 47[e], date d'avril 2009. Il s'agit de *$2^{42\,643\,801} - 1$*, et donc le 47[e] (chronologiquement) nombre parfait associé (suivant la proposition IX.36) est : *$2^{42\,643\,800}(2^{42\,643\,801} - 1)$*. On doit parler d'un ordre chronologique, car on ne sait pas s'il n'y a pas de nombres intermédiaires. Comme on le voit, on a affaire très rapidement à des nombres colossaux (plusieurs millions de chiffres).

[84] Platon, *Lois* [= PLA4], 919e6-920a3 (trad. Brisson-Pradeau).

[85] Voir par exemple la quasi-totalité du *Sophiste* ; *Pol.* 148c ; *Phaedr.* 266a et c ; Vuillemin 2001, p. 23, propose même d'établir « une correspondance bi-univoque entre l'interprétation du tableau (platonicien) de Nicomaque et l'arbre dichotomique de Platon ».

[86] Arist. *P.A.* I, 643a20-25.

[87] *Metaph.* M, 7, 1082b3.

[88] *Metaph.* D, 29, 1025a1 ; M, 7, 1081a14-25, b10-30 ; M, 8, 1083a13-14, b30-36 ; M, 9, 1085b5-10 ; N, 3, 1091a4-12. C'est à une conclusion identique qu'est amené Vuillemin 2001 dans son étude de la dyade indéfinie (p. 16-17). Et, partant du *Phédon*, 105e, et des textes de la *Métaphysique* d'Aristote, c'est encore à quoi aboutit Wedberg 1955, p. 137-138 : « une conjecture plausible est que [les nombres] doivent être engendrés par doublement, par la multiplication par *2*. Si cette interprétation est correcte, Platon pense ici à la série des nombres comme engendrée à partir du nombre initial *2* par (i) l'opération *n + 1* appliquée à tout nombre pair, et (ii) l'opération *2 × n* appliquée à un nombre quelconque *n*. *I.e.*, la série des nombres est construite de la manière suivante : *2, 2 + 1, 2 × 2, 2 × 2 +1, 2 × (2 + 1), …* » ('*A plausible guess, however, is that [the numbers] are to be generated through doubling, through multiplication by* 2. *If this interpretation is correct, Plato is here thinking of the number series as generated from the initial number* 2 *by (i) the operation* n+1 *being applied to any given* n, *and (ii) the operation* 2 × n *being apply to any* n. *I.e., the number series is construed as follows:* 2, 2 + 1, 2 × 2, 2 × 2 +1, 2 × (2 + 1), …'.

[89] *Metaph.* M, 7, 1081a14-15.

*le nombre pair provient ou bien d'une duplication de la Dyade à partir de l'Un ou bien d'une duplication des nombres impairs <engendrés par l'application de l'Un au nombre pair>*[90].

Ce qui n'est rien d'autre que la propriété pour tout pair de s'écrire sous la forme d'une certaine puissance de *2*, multipliée par un impair, *i.e.* la décomposition pair/impair[91].

La démonstration de cette propriété est beaucoup plus simple que celle, par l'impossible, de la proposition VII.22. Elle se ramène d'ailleurs à des résultats très primitifs, induits par la dichotomie itérée, qui sont parmi les plus anciens connus en mathématique, puisque remontant aux Babyloniens[92] et à l'Égypte antique.

Non seulement la série des puissances de *2* y était communément utilisée pour la multiplication entière, mais aussi celle de ses inverses. Ainsi, selon une tradition en égyptologie, cette série jusqu'à l'ordre *6* (i.e. *1/64*) prenait le nom de « fractions œil d'Horus » (*oudjat* en égyptien)[93]. Et lorsqu'il repère la réforme métrologique introduite par Solon à partir des puissances de *2*, M. Caveing y voit une « application judicieuse des quantièmes égyptiens "œil d'Horus" ». Il considère également qu'au VI$^e$ siècle BCE, l'usage de la duplication ou de la division par *2* « était usuel et dominant, sinon exclusif » chez les mathématiciens grecs[94].

On remarque en effet que la multiplication par des puissances de *2* est particulièrement simple dès lors que l'on décompose le multiplicateur en somme de telles puissances, ce qui est toujours possible[95]. Ainsi, *7* s'écrit : *7 = $2^2$ + 2 + 1* et *10 = $2^3$ + 2*.

Suivant le papyrus de Rhind, cette écriture, comme somme de puissances de *2*, était au fondement de la technique multiplicative égyptienne. Pour multiplier deux nombres, on décompose le multiplicateur sous la forme d'une somme de telles puissances, puis on procède à la multiplication par chaque élément de cette somme, et enfin on additionne l'ensemble[96].

---

[90] *Ibid*. M, 8, 1084a4-7 (trad. Tricot) [=ARI1].

[91] *Cf.* aussi J. Tricot ([= ARI1]), II, p. 240 note 2 et Caveing 1998, II, p. 209. Knorr propose une origine pythagoricienne plus géométrique, qui serait également la source des recherches sur l'incommensurabilité entreprises par les mathématiciens grecs anciens. On considère un triangle rectangle isocèle de côtés tous entiers et pairs. Il établit tout d'abord un théorème (numéroté par lui V.14) montrant que si l'hypoténuse d'un tel triangle est pair, alors il en est de même des autres côtés. On construit alors le triangle obtenu en joignant les milieux de tous les côtés. C'est encore un triangle rectangle isocèle (*cf.* figure ci-dessous) à côtés entiers, et si l'hypoténuse est paire, on peut itérer.

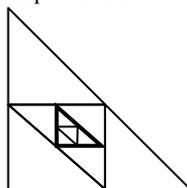

La division peut être poursuivie tant que l'hypoténuse est paire. Les longueurs des côtés du deuxième triangle sont égales à celle du premier divisées par *2*, celles du troisième par $2^2 = 4$, du quatrième par $2^3 = 8$, du cinquième par $2^4 = 16$, … « Puisque une régression indéfinie contredit sa nature d'être un entier » ('*since an indefinite regress contradicts the integral nature of [the hypotenuse]*'), ce processus doit s'arrêter. Ce qui permet de réduire la recherche au seul cas où l'hypoténuse correspond à un impair, et de là, établir une contradiction (Knorr 1975, p. 179-180). Cela revient donc à diviser l'hypoténuse par *2* autant de fois que possible, pour finalement obtenir un impair. Ici encore, la division finie par *2* est considérée comme une évidence pour la mathématique grecque du temps de Pythagore.

[92] Knorr 1975, p. 23 note 53, souligne que l'on trouve des exemples de dichotomies itérées sur les tablettes mésopotamiennes antérieures à 1600 BCE.

[93] *Cf.* par exemple Allen 2000, p. 101. On a toutefois récemment montré la fragilité des témoignages sur lesquels repose cette interprétation (Jim Ritter, « Closing the Eye of Horus … », dans A. Imhausen et J. Steele (éd.), *Under One Sky : astronomy and mathematics in the ancient near East*, Münster, 2003, p. 298-323). Ce qui importe toutefois ici est l'existence d'un récit de la mythologie égyptienne, où le dieu Horus eut l'œil arraché dans un combat avec son frère Seth, éclaté en une suite de morceaux, chacun moitié moindre du précédent, puis restauré par Thot.

[94] Caveing 1998, II, p. 201.

[95] Pour M. Caveing, le « théorème fondamental », quoique implicite, de l'arithmétique égyptienne établit que tout entier naturel est une somme (finie) et unique d'éléments de la suite des puissances de *2* (auxquelles s'ajoute l'unité pour les nombres impairs) (Caveing 1998, II, p. 202).

[96] Granger 1976, p. 27, remarque l'importance d'une technique « où la multiplication d'un entier est ramenée à l'addition des produits de celui-ci par celles des puissances de 2 successives en lesquelles on peu toujours décomposer la

En reprenant l'exemple donné par Heath[97], pour multiplier, dans notre système décimal, un nombre par *13*, nous utilisons sa décomposition en base *10* (i.e. *13 = 10 + 3*). Nous multiplions alors ce nombre par *10*, puis (ce même nombre) par *3*, et nous additionnons les deux résultats ainsi obtenus.

Dans le système égyptien, pour effectuer cette multiplication par *13* (= *8 + 4 + 1* = $2^3 + 2^2 + 1$), on multiplie le nombre par *8*, puis (ce même nombre) par *4*, puis enfin par *1*, et on additionne alors les trois résultats[98]. Par exemple, en écriture moderne, la multiplication de *15* par *13* aurait la forme suivante :

*15 × 13 = 15 × (8 + 4 + 1) = 15 × 8 + 15 × 4 + 15 × 1 = 120 + 60 + 15 = 195*.

Certes, la manière usuelle grecque de multiplier était différente de celle des Égyptiens. Néanmoins, suivant un scholiaste du *Charmide* de Platon, les anciens Grecs en utilisaient également une autre qu'il nomme « méthode égyptienne » par opposition à la « méthode grecque ». Il aurait d'ailleurs été tout à fait extraordinaire qu'une technique calculatoire aussi simple et ancienne ait complètement échappé au Grecs anciens dont, dès la haute antiquité, les colonies africaines jouxtaient l'Égypte. Les plus anciens textes grecs notent d'ailleurs l'usage des Athéniens cultivés de séjourner en Égypte, afin d'y retrouver les sources du savoir. Ainsi Platon y situe, entre autre, l'origine mythique des mathématiques (*Phaedr.* 274c-275d). Et dans les *Lois*, c'est l'éducation mathématique des enfants égyptiens qu'il propose comme modèle pour celle des jeunes Grecs (VII, 819a-d).

Maurice Caveing, comparant les constructions égyptiennes à celle donnée par Aristote, remarque leur similitude[99]. La première consiste en une décomposition complète (*i.e.* l'écriture du nombre en « base *2* »), la seconde à diviser autant de fois que possible par *2*. Celle-là implique celle-ci très simplement, il suffit de considérer la plus petite puissance de *2* qui intervient dans la décomposition du nombre. Inversement, par itération, on passe de la seconde à la première.

Sous sa forme négative, la décomposition pair/impair énonce que, pour tout nombre donné, la division successive par *2* ne peut se poursuivre indéfiniment. Elle est utilisée sous cette forme, généralement de manière implicite, dans les démonstrations élaborées par les commentateurs pour déduire l'irrationalité de la racine de *2* par le pair et l'impair. C'est elle qui permet de distinguer les nombres des grandeurs « continues », celles-ci étant infiniment divisibles en deux parties (égales) contrairement à ceux-là[100]. Ainsi, M. Caveing donne cette propriété, mais elle apparaît comme une conséquence explicative de l'impossibilité d'une descente infinie pour les entiers, plutôt qu'un résultat obtenu pour lui-même[101].

Un autre texte des *Premiers analytiques* rapproche encore les arguments de Zénon sur la dichotomie itérée (donc les puissances inverses de *2*) et l'irrationalité :

*si par exemple, voulant prouver l'incommensurabilité de la diagonale, on cherchait à démontrer l'argument de Zénon sur l'impossibilité du mouvement, et qu'en vue de cette*

---

multiplication ». Puis il relie cette manière d'opérer à Platon et à Aristote : « Platon, on le sait, donne un sens métaphysique à la genèse du nombre par duplication, et Aristote, discutant Platon (*Métaphysique*, M, 1081b), indique à ce propos une distinction particulièrement intéressante entre la "dyade première" et la "dyade indéterminée" ».

[97] Heath 1981, p. 52-53.

[98] En termes modernes, pour le multiplicateur, on utilise une représentation binaire des entiers. Cela suppose la connaissance du fait que tout entier (naturel) peut s'écrire, de manière unique, comme somme de puissances de *2* (Caveing 1998, II, p. 201-209). Pour d'autres exemples, *cf.* Ifrah 1994, I, p. 423-427.

[99] Caveing 1998, II, p. 211.

[100] Lorsque Penrose recherche une preuve de l'irrationalité de $\sqrt{2}$ aussi simple que possible pour le lecteur contemporain (cf. *supra* note 73), il récuse la démonstration usuelle (i.e. *via* la proposition VII.22 des *Éléments*) au profit de cette propriété : l'impossibilité d'obtenir une suite d'entiers (distincts) de cette forme (*i.e.* par divisions successives par *2*). De l'hypothèse $\sqrt{2} = a/b$ où *a* et *b* sont des entiers (positifs non nuls), il argumente de la manière suivante : « Nous pouvons donc reproduire le même raisonnement encore et encore, pour obtenir une succession sans fin d'équations $a^2 = 2b$, $b = 2c$, $c = 2d$, $d = 2e$, … (…) Or toute série décroissante d'entiers positifs doit avoir une fin tandis que celle-ci n'en a pas. Nous aboutissons donc à une contradiction' (Penrose 2007, § 3.1, p. 48).

[101] Caveing 1998, III, p. 130.

> *proposition on procédât par la réduction à l'absurde : absolument d'aucune façon, la conclusion fausse n'est en connexion avec l'énonciation du début.* (*An. pr*. II, 17, 65b16-21, trad. Tricot).

Certes Aristote donne ici un contre-exemple à un véritable syllogisme, consistant à vouloir prouver « l'incommensurabilité de la diagonale » par réduction à « l'impossible » de l'argumentation de Zénon sur le mouvement. Mais en les rapprochant, il signifie aussi qu'il existe quelque chose de commun entre eux. La reconstruction qu'en propose Caveing se fonde précisément sur la division infinie par $2$[102]. Une fois encore, c'est la suite des puissances de *2* (ou de ses inverses) qui apparaît au cœur de questions très anciennes.

---

[102] Caveing 1982, p. 119. Knorr 1975, p. 27, formule essentiellement le même point de vue sur cette question.

## II. Une propriété essentielle de la décomposition pair/impair

(a) Cette propriété concerne le produit de deux nombres :
**le degré de parité d'un produit est la somme des degrés de parité de chacun des termes de ce produit** *i.e.*
$m = 2^k u$ **(avec *u* impair) et** $n = 2^h v$ **(*v* impair) ;**
alors
$m \times n = 2^{k+h} w$ **(où *w* est impair**, et en fait égal au produit de *u* et *v* i.e. $w = u \times v$**)**.
Nous appellerons cette propriété l'**additivité des degrés de parité**[103].
En particulier, suivant cette propriété, le degré de parité de $m^2 = m \times m$ (*i.e.* le nombre de fois tel que $m^2$ est divisible par *2*) est le double du degré de parité de *m*. On a donc la propriété suivante :
**Le degré de parité d'un carré est toujours pair**, ce que nous appellerons la **propriété des degrés de parité des carrés**.
De même, **le degré de parité d'un cube est toujours un multiple de 3**, ce que nous appellerons la **propriété des degrés de parité des cubes**.

(b) De manière plus formelle, on peut voir la propriété d'additivité des degrés de parité de la façon suivante :
1. Tout d'abord, si un nombre *n* est divisible *h* fois par *2*, alors *2n* est divisible *h+1* fois.
En termes moderne, si $n = 2^h u$ (où *u* est impair), alors $2n = 2^{h+1} u$.
Exemples :
(a) $8 = 2 \times 2 \times 2 = 2^3$ (*h = 3*) d'où $2 \times 8 = 2 \times 2 \times 2 \times 2 = 2^4$ (et *4 = 3+1 = h+1*) ;
(b) $14 = 2 \times 7$ (*h = 1*), d'où $2 \times 14 = 2 \times 2 \times 7 = 2^2 \times 7$ (et *2 = 1+1 = h+1*).
2. De même, $2^k n$ est divisible *h+k* fois *i.e.* $2^k n = 2^{h+k} u$.
Exemples :
soit $8 = 2 \times 2 \times 2 = 2^3$ (*k = 3*) :
(a) pour $n = 8 = 2^3$ (*h = 3*), on a :
$8 \times n = 8 \times 8 = (2 \times 2 \times 2) \times (2 \times 2 \times 2) = 2^6$ (et *6 = 3+3 = h+k*) ;
(b) pour $n = 14 = 2 \times 7$ (*h = 1*), on a :
$8 \times n = 8 \times 14 = 2^3 \times 14 = (2 \times 2 \times 2) \times (2 \times 7) = 2^4 \times 7$ (et *4 = 3+1 = h+k*).
3. Et donc, si *m* est un entier divisible *k* fois par *2*, et *n* un entier divisible *h* fois par *2*, leur produit $m \times n$ est divisible *k+h* fois par *2*.
Exemples :
(a) Soit $n = 4 = 2 \times 2 = 2^2$ (*h = 2*) et $m = 12 = 4 \times 3 = 2^2 \times 3$ (donc *k = 2*).
On a alors :
$12 \times 4 = (2^2 \times 3) \times 2^2 = 2^4 \times 3 = 2^{h+k} \times 3$.
(b) Soit $n = 6 = 2 \times 3$ (*h* = 1) et $m = 12 = 4 \times 3 = 2^2 \times 3$ (donc *k = 2*).
On a alors :
$m \times n = 12 \times 6 = (2^2 \times 3) \times (2 \times 3) = 2^3 \times 9 = 2^{h+k} \times 9$.
(c) Soit $m = 6 = 2 \times 3$ (*h = 1*). On a alors :
$m^2 = m \times m = 6 \times 6 = (2 \times 3) \times (2 \times 3) = 2^2 \times 9$,
et le degré de parité de $6^2 = 36$ est *2* qui est bien le double du degré de parité de *6*.
(d) Soit $n = 2$ (*h = 1*). On a alors :

---

[103] Le cas où l'un des deux nombres *m* ou *n* est impair est particulièrement simple. S'ils ne sont pas de même parité, le degré de parité de leur produit est celui du nombre pair. Cela résulte immédiatement de ce que le produit de deux nombres impairs est impair (par exemple si $n = 2^h u$ (*u* impair) et *m* est impair, $n \times m = 2^h (u \times m)$, et son degré de parité est *h*). Enfin si les deux nombres sont impairs, donc non divisibles par *2*, leur produit est impair, c'est-à-dire toujours non divisible par *2*.

$m^3 = m \times m \times m = 2^3 = 8$,

et le degré de parité de *8* est *3*, triple du degré de parité de *2*.

Cette propriété résulte pratiquement de la définition du degré de parité. Elle est d'ailleurs le fondement de la multiplication utilisée par les Égyptiens de l'antiquité, et connue par les Grecs anciens, au moins depuis les pythagoriciens anciens (sans doute parce que rien ne nous est parvenu, ne serait-ce qu'indirectement, des mathématiques grecques antérieures).

# III. La démonstration de l'incommensurabilité de la diagonale au côté du carré (*i.e.* de l'irrationalité de √2)

Comme dans les démonstrations standard, on va procéder au moyen d'un raisonnement par l'impossible.

Il s'agit de montrer que l'existence d'entiers *m* et *n* de rapport égal à √2, conduit à une impossibilité.

Supposons donc que l'on ait :
$\sqrt{2} = m/n$.

En élevant au carré, on obtient :
$2 = (m/n)^2 = m^2/n^2$

ou encore :
$2n^2 = m^2$ (égalité 1).

Suivant l'additivité des degrés de parité (cf. *supra* II.a), le degré de parité de *2n²* est égal à celui de *n²* augmenté d'une unité. Et suivant la propriété des degrés de parité des carrés (*ib.*), les degrés de parité de *n²* et de *m²* sont pairs.

Les membres de gauche et de droite de l'égalité représentant le même nombre, les degrés de parité de droite et de gauche (*i.e.* le nombre maximal de divisions par *2* de ce nombre) sont donc identiques.

**On a donc : un nombre pair augmenté d'une unité est égal à un nombre pair.**

L'unité, différence de deux pairs, est donc paire[104].

Par définition, un impair est la somme d'un pair et d'une unité[105]. Celle-ci étant paire, tout impair est alors somme de deux pairs. La somme de deux pairs étant paire, on obtient que tout impair est pair, et donc l'effondrement de l'impair dans le pair, ce qui est bien la contradiction, indiquée par trois fois, dans les *Premiers Analytiques* d'Aristote.

De l'égalité 1, on pourrait certes aboutir à une égalité entre *un* nombre impair et *un* nombre pair. Mais contrairement aux démonstrations usuelles, qui raisonnaient abstraitement sur les nombres, l'unité apparaît ici explicitement. C'est ce qui permet de conclure non pas (simplement) qu'un nombre (indéterminé) possède le *double attribut* d'être impair et pair, mais que l'unité est paire[106].

Cette démonstration n'utilise pas la proposition VII.22 des *Éléments*, qui se trouve au cœur de la preuve générale d'irrationalité des entiers non carrés. Une fois celle-ci obtenue, l'impossibilité à laquelle aboutissait le raisonnement se déplaçait de la parité à la primalité[107], et le résultat particulier concernant *2* ainsi que la méthode fondée sur la parité devaient alors disparaître. Cela explique que, dès l'antiquité, mais quelques siècles après Aristote, les commentateurs se soient évertués à le réconcilier avec Euclide, en proposant d'autres démonstrations fondées sur cette proposition. Du point de vue mathématique, la tendance est, en effet, à la généralisation et à l'oubli des étapes historiques intermédiaires.

---

[104] L'unité serait alors divisible par *2*, ce qui ne ferait pas rire les mathématiciens de la *République* dont la pire crainte serait réalisée, l'unité n'apparaissant « plus comme un, mais comme des parties multiples » (μὴ ἓν ἀλλὰ πολλὰ μόρια) (VII, 525e5). Pour les Grecs anciens, l'unité était précisément conçue, par définition, comme indivisible. D'autre part, puisque *1* est le plus petit impair et *2* le plus petit pair, on aurait *1 = 2*, ce qui paraîtrait absurde même au dormeur de Platon (cf. *Theaet.* 190c ; et aussi *supra*, première partie, chap. V).

[105] C'est la définition 7 du livre VII des *Éléments* d'Euclide : « Un impair est un nombre (…) qui diffère d'une unité d'un nombre pair. » Mais on la trouve déjà dans le *Phédon* de Platon : « c'est la *monade* qui *rendra impair* tout ce en quoi elle s'insère. » ( Dixsaut 1991.105c). Th. Heath donne la définition de l'impair dérivée « de la conception populaire » des Grecs anciens (ἐκ δεμώδους ὑπολήψεως), telle que rapportée par Nicomaque : « est pair ce qui est capable d'être divisé en deux parts égales sans une unité tombant au milieu, et est impair ce qui ne peut être divisé en deux parties égales à cause de ladite intervention des unités » ('*that is even which is capable of being divided into two equal parts without an unit falling in the middle, and that is odd which cannot be divided into two equal parts because of the foresaid intervention of the units*'(Heath 1956, II, p. 281).

[106] Maurice Caveing aboutit lui aussi à cette formule. Mais c'est pour en déduire que *p* et *q* ne peuvent être simultanément des puissances de *2*. Il l'obtient au moyen d'un nouveau raisonnement par l'impossible, aboutissant encore à ce qu'un nombre pair deviendrait alors égal à un nombre impair (Caveing 1998, III, p. 130).

[107] Cf. *supra*, première partie, chap. V, en particulier notes 47 et 49.

La diversité des preuves « usuelles » supposées rendre compte des énoncés aristotéliciens sur l'incommensurabilité résulte de la manière dont on les lit : si l'on s'attache aux termes, on essaiera d'introduire, de manière plus ou moins raffinée, la notion de parité. Si l'on s'intéresse plutôt à la cohérence mathématique, on sera amené à s'en débarrasser.

Dans la construction que nous avons proposée, au contraire, la parité joue un rôle crucial, *via* la décomposition pair/impair. Dans la démonstration, elle intervient de manière très spéciale, concluant à ce que l'unité devient paire. D'où s'ensuit l'effondrement de la notion même de parité, l'impair étant défini comme l'addition d'un pair et de l'unité (cf. *supra* note 103).

On peut dès lors interpréter littéralement ces textes des *Premiers analytiques*. On conclut en effet, non pas à ce qu'*un* impair soit égal à un pair, ou encore à ce qu'*un* nombre soit à la fois pair et impair, ce à quoi aboutissent les autres démonstrations, mais à ce que tout impair est pair. Ou encore, en termes aristotéliciens, que (le genre) impair est (le genre) pair (τὰ περιττὰ ἴσα εἶναι τοῖς ἀρτίοις), les neutres pluriels (περιττα, ἀρτίοις) désignant la classe générale définie par le nom[108].

---

[108] Pour revenir à la discussion sur la symétrie ou non de la relation d'égalité dans les *Premiers analytiques* (cf. *supra*, première partie, chap. V, en particulier notes 60 et 64), dès que l'unité est paire, on a confusion entre impairs et pairs. En effet, tout pair s'écrit alors comme un pair (le nombre qui le précède) augmenté de l'unité. Ainsi 2 est égal à la somme de *1* (qui est pair) et de *1* ; de même, *4 = 3 + 1* où *3* est un pair (puisque tout impair est pair) augmenté de *1*, et *4* est donc impair. Autrement dit, dès que les impairs deviennent pairs, tous les nombres deviennent à la fois pairs et impairs, ce qui ôte toute signification à ces termes. On a toutefois besoin d'une propriété non utilisée dans notre démonstration (ainsi *2 = 1 + 1*).

# IV. Conséquences

**(a) Le résultat reste vrai pour tout nombre pair dont le degré de parité est impair**

Il suffit de répéter la démonstration précédente.

En effet, soit $x = 2^i s$ un entier, où $i$ et $s$ sont impairs (décomposition pair/impair de $x$).

Supposons par l'absurde qu'il existe des entiers $m$ et $n$ dont le rapport est égal à $\sqrt{x}$.

En élevant au carré, on obtient :

$x = (m/n)^2 = m^2/n^2$

ou encore :

$x \times n^2 = m^2$ ⟶ (égalité 1')

Cette égalité signifie que les membres de gauche ($x \times n^2$) et de droite ($m^2$) représentent le même nombre.

Suivant l'additivité des degrés de parité (cf. *supra* II.a), le degré de parité du membre de gauche est égal à la somme de ceux de $x$ et de $n^2$. Les degrés de parité de $n^2$ et $m^2$ étant pairs (*ibid.*), la somme d'un nombre impair (l'entier $i$) et d'un pair (le degré de parité de $n^2$) est égal à un nombre pair (le degré de parité de $m^2$).

La somme d'un impair et d'un pair étant impaire, on a une égalité entre un pair et un impair, d'où une impossibilité[109].

**Application**

Les racines carrées de tous les nombres pairs tels *2* (ce qu'on sait déjà), *6, 8, 10, 14* (mais ni *4*, ni *12*, ni *16*) sont irrationnelles[110].

**(b) L'irrationalité ramenée au seul cas des nombres impairs**

Nous utiliserons ici une propriété très simple sur les carrés :

si on multiplie le côté d'un carré par un entier, son rapport au premier carré est le carré de cet entier.

Ainsi, en doublant le côté d'un carré, on obtient un carré *4* fois plus grand, comme Socrate le montre à la surprise du jeune serviteur du *Ménon* de Platon (82e-83b), qui pensait ainsi doubler le carré initial.

Puisque $\sqrt{t}$ est le côté du carré de surface $t$, en écriture moderne, cela donne :

$\sqrt{a^2 t} = a\sqrt{t}$ ⟶ (égalité 2)

où $a$ et $t$ sont des entiers.

Considérons maintenant le cas d'un entier $y$ de degré de parité pair, *i.e.* sa décomposition pair/impair donnée par :

$y = 2^{2i} s$, où $s$ est impair,

et, d'après l'additivité des degrés de parité (cf. *supra* II.a), on a :

$(2^i)^2 = 2^i \times 2^i = 2^{2i}$.

**1.** La racine carrée de $y$ est rationnelle signifie qu'il existe des entiers $p$ et $q$ tels que :

$\sqrt{y} = p/q$ ou encore (puisque $y = 2^{2i} s = (2^i)^2 s$),

---

[109] À la différence du paragraphe précédent, on conclut sur cette contradiction, plutôt que sur l'unité comme paire. On a en effet affaire, cette fois, à la somme d'un impair et d'un pair, et non pas de l'unité et d'un entier pair. Parce que ce résultat englobe, en tant que cas particulier, celui $x = 2$, il pourrait avoir fait « oublier » aux mathématiciens postérieurs la conclusion impossible de la démonstration primitive (l'unité est paire), au profit de l'égalité d'un pair et d'un impair que l'on retrouve dans les démonstrations standard.

[110] Le sens de la parenthèse est quelque peu ambigu. On ne veut pas signifier que les racines carrées de ces trois entiers ne sont pas irrationnelles, ce qui est d'ailleurs faux, mais que cette négation ne résulte pas de ce qui précède.

$\sqrt{y} = \sqrt{2^{2i} s} = \sqrt{(2^i)^2 s} =$ (suivant l'égalité 2) $2^i \sqrt{s} = p/q$

i.e. $\sqrt{s}$ est rationnel.

La rationalité de $\sqrt{y}$ étant équivalente à celle de $\sqrt{s}$, leur irrationalité l'est aussi.

Autrement dit, pour un entier pair $z = 2^j w$ (où $w$ est impair) deux cas sont possibles :

– ou bien son degré de parité $j$ est impair, et alors sa racine carrée est irrationnelle ;

– ou bien $j$ est pair, et dans ce cas, $\sqrt{z}$ est rationnel (ou irrationnel) selon que $\sqrt{w}$ l'est (ou ne l'est pas). Ainsi, la rationalité (ou pas) de racine carrée de *12 = 4 × 3 = $2^2$ × 3* se ramène à celle de *3*.

Pour résoudre le problème de la rationalité/irrationalité des racines carrées d'entiers, il suffit donc de résoudre celui des seuls nombres impairs.

**2.** En outre, la racine carrée de *w* (entier impair) est rationnelle équivaut à *w* est le carré d'un nombre rationnel, *i.e.* il existe des entiers *m* et *n* tels que :

*w = (m/n)² = m²/n²*

ou encore :

*w × n² = m²*                                                                                      (égalité 3).

Si on note $m = 2^h u$ et $n = 2^k v$ (*u* et *v* impairs), les décompositions pair/impair de *m* et *n*, d'après l'additivité des degrés de parité (cf. s*upra* II.a), celui de $n^2$ est *2h*, et celui de $m^2$ est *2k*, i.e.

$m^2 = 2^{2h} t$ et $n^2 = 2^{2k} r$ où $t = u^2$ et $r = v^2$ sont impairs.

L'entier *w* étant impair, le degré de parité de $w \times n^2$ est égal à celui de $n^2$ i.e. à *2h* (cf. *supra* II, note 101). D'après l'égalité 3, on a donc :

*2h = 2k* i.e. *h = k*,

d'où :

$w = m^2/n^2 = (2^{2h} t)/(2^{2k} r) = t/r = u^2/v^2 = (u/v)^2$,

et *w* est le carré d'un rapport de deux impairs.

Ainsi, d'après les numéros 1 et 2 ci-dessus, *la question portant sur les racines carrées est entièrement ramenée au cadre des entiers impairs.*

*En effet, il s'agit d'étudier uniquement les racines carrées des impairs, mais pour être rationnelles, celles-ci doivent en outre être quotient de deux entiers impairs.* Par exemple, la rationalité de racine carrée de *12* se ramène à celle de *3* ; mais, *en outre*, la racine carrée de *3*, pour être rationnelle, doit être le quotient de deux impairs.

Dans le cadre de cette démonstration, il est inutile de recourir à la proposition VII.22 des *Éléments*, selon laquelle tout rationnel peut être représenté comme rapport de deux entiers *premiers* entre eux, nombres minimaux parmi tous les rapports donnant ce rationnel. Cette proposition n'est aucunement triviale, et pose un problème chronologique difficile (cf. *supra* note 49). La faire intervenir ajoute un degré d'incertitude très fort pour dater la connaissance de l'irrationalité de $\sqrt{2}$, et par conséquent celle des autres racines irrationnelles d'entiers. Inversement, sa datation influencera celle de la proposition VII.22, voire de la totalité du livre VII, avec, une fois de plus, un risque de circularité évident.

Suivant notre démonstration, au contraire, ces deux questions sont indépendantes.

**(c) Les cas des racines cubiques**

Cette méthode du raffinement du pair et de l'impair s'étend facilement au cas des racines cubiques.

Au paragraphe précédent, on a donné une condition nécessaire pour que la racine carrée d'un entier soit rationnelle (*i.e.* que le côté d'un carré soit commensurable à l'unité) : le degré de parité de cet entier doit être pair. Ou de manière négative (cf. supra, note 71), si le degré de parité d'un entier est impair, sa racine carrée est irrationnelle.

On va donner ici une condition nécessaire pour la rationalité de la racine cubique d'un entier (*i.e.* pour que le côté d'un cube de volume entier soit commensurable à l'unité).

Soit *x* un entier pair. Supposons qu'il existe des entiers *m* et *n* dont le rapport est égal à la racine cubique de *x*, i.e.
$\sqrt[3]{x} = m/n$.

En élevant au cube, des deux côtés de l'égalité, on obtient :
$x = (m/n)^3 = m^3/n^3$
ou encore :
$x \times n^3 = m^3$ (égalité 1'').

D'après la propriété du degré de parité des cubes (cf. *supra* II.a), on a :
le degré de parité de $n^3$ et $m^3$ est un multiple de *3* (propriété 1)
et d'après la propriété d'additivité des degrés de parité (*cf.* toujours II.a) :
le degré de parité de $x \times n^3$ est la somme du degré de *x* et de celui de $n^3$ (propriété 2).

D'après l'égalité 1'', les degrés de parité de $x \times n^3$ et de $m^3$ sont identiques.

D'après la propriété 2, le degré de parité de *x* est donc la différence des degrés de parité de $m^3$ et $n^3$.

Et donc, d'après la propriété 1, le degré de parité de *x* est la différence de deux multiples de *3*, c'est donc un multiple de *3*.

Ainsi une condition nécessaire pour que la racine cubique d'un entier soit rationnelle est que son degré de parité soit un multiple de *3*. Ce qui s'écrit encore négativement (cf. supra, note 71) :

**Propriété des racines cubiques**

Si le degré de parité d'un entier n'est pas un multiple de *3*, sa racine cubique est irrationnelle.

**Application**

Les racines cubiques de *2, 4, 6, 10, 12, 14, 16, 18, 20, 22* ne peuvent être rationnelles. Pour celle de *24 = 8 × 3 = 2³ × 3*, on est ramené au problème de la racine cubique de *3*[111]. Ainsi, cette fois encore, le problème se ramène à la seule question des impairs[112].

---

[111] On utilise l'analogue pour les cubes de la propriété (2) des carrés (cf. *supra*, b) : le rapport des volumes de deux cubes est le cube du rapport de leurs côtés, ce qui, pour *a* entier, donne en écriture moderne : $\sqrt[3]{a^3 x} = a\sqrt[3]{x}$.

[112] On obtient ainsi un parallélisme complet avec le cas des racines carrées, comme indiqué dans le *Théétète* (148b ; pour une analyse de ce passage, voir notre article à paraître dans un prochain numéro.

# V. Sur l'adéquate généralité de la démonstration.

Nous avons souligné que l'utilisation de la proposition VII.22 dans les démonstrations standard posait un problème vis-à-vis de la théorie aristotélicienne de la science. Elle est beaucoup trop générale et (*via* la proposition VII.27 qui en est en quelque sorte le pendant) elle donne la solution pour tous les entiers non carrés (d'entiers). Ainsi s'expose-t-elle à une double critique. Épistémologique, en tant que contrevenant à l'adéquation exigée par Aristote entre le moyen (la preuve) et la fin (le résultat obtenu) ; chronologique, en tant que les témoignages textuels imposent une période suffisamment longue entre le résultat particulier de l'irrationalité de racine carrée de *2*, et le cas général.

Reste à voir si notre démonstration ne s'expose pas, elle-même, à cette critique.

De fait, il ne paraît pas très difficile de passer de la preuve particulière donnée au chapitre IV, à la preuve générale[113].

Soit donc *x* un entier *quelconque*. On reprend la démonstration du chapitre III, en remplaçant *2* par *x*.

Tout entier *n* peut s'écrire sous la forme d'un produit $n = x^h u$ où *h* est le nombre maximal de divisions possibles de *n* par *x*, et *u* est un entier non divisible par *x*.

De la définition de *h* suit immédiatement son unicité (et donc celle de *u*), ainsi que les propriétés d'additivité des degrés (cf. *supra* chap. II).

On suppose alors que la racine carrée de *x* est rationnelle, *i.e.* elle s'écrit :
$\sqrt{x} = m/n$ avec $n = x^h u$ et $m = x^k v$ où *u* et *v* sont des entiers non divisibles par *x*,
d'où :
$x = m^2/n^2$ et donc $x \times n^2 = m^2$,
et (car $n = x^h u$ et $m = x^k v$), on obtient finalement :
$x^{2h+1} t = x^{2k} s$ avec $t = u^2$ et $s = v^2$ (résultat 1).

Pour conclure par unicité, comme dans le cas *x = 2*, il suffirait alors que *t* et *s* ne soient pas divisibles par *x*. Mais on sait seulement que *u* et *v* ne le sont pas, il s'agit donc de passer de ceux-ci à ceux-là.

Pour *x = 2*, cela est évident, car les entiers et leurs carrés ont même parité. Il n'en est pas de même dans la situation générale qui, au contraire, nécessite une arithmétique suffisamment développée, tout à fait incompatible avec l'époque que nous avons considérée[114].

Cela est d'ailleurs faux en général : qu'un entier ne soit pas divisible par un autre ne signifie pas qu'il en soit de même de son carré[115].

---

[113] D'après ce qui précède, la question se pose pour les seuls entiers impairs. Mais cela ne joue pas de rôle dans l'argumentation qui suit, et nous considérons tous les entiers y compris les carrés d'entiers. Comme dans la leçon de Théodore, dans le *Théétète*, il s'agit de distinguer entre la méthode de démonstration qui s'applique à tous les entiers, et le résultat qu'on veut obtenir qui, lui, n'est vrai que pour certains d'entre eux (ceux qui ne sont pas égaux au carré d'un entier).

[114] Ainsi, Th. Heath, qui conteste la thèse selon laquelle les pythagoriciens anciens n'auraient pas eu une (certaine) théorie des proportions, considère qu'elle devait être très rudimentaire en comparaison de celle qu'on trouve aux livres VII et VIII des *Éléments* (*cf.* Heath 1981, p. 211-212).

[115] Ainsi *4* n'est pas divisible par *8*, mais son carré $4^2 = 16$ l'est ; *18* n'est pas divisible par *12*, pourtant son carré l'est (car $18^2 = 324 = 12 \times 27$) et *15* n'est pas divisible par *9*, mais son carré $15^2 = 225 = 9 \times 25$ l'est. Il faut donc élaborer une théorie pour résoudre cette question, ce qui est fait à la fin du livre VII des *Éléments*, *via* en particulier les propositions VII.26-32. C'est essentiellement la théorie de la décomposition des entiers en facteurs premiers qui est alors élaborée. Mais de celle-ci découle immédiatement le résultat général pour les racines carrées d'entiers quelconques.

On a vu dans la preuve donnée par Alexandre, que celui-ci utilisait des carrés au lieu des entiers de la démonstration standard (*cf. infra*, première partie, §II.2). Cela irait dans le sens d'une antériorité de la preuve d'Alexandre sur celle qu'on trouve à la fin du livre X des *Éléments*.

Et de fait, pour conclure le raisonnement par l'absurde prouvant l'irrationalité de √*x*, il faut montrer l'impossibilité que le résultat 1 soit vérifié (quels que soient les entiers *k, h, u* et *s*). Or celui-ci peut se réécrire de la manière suivante :
- Si *2h+1* est plus grand que *2k* : $x^{2h+1-2k} = s/t = v^2/u^2$ i.e. $\sqrt{x^{2h+1-2k}} = v/u$
- Si *2k* est plus grand que *2h+1* : $x^{2k-2h-1} = t/s = u^2/v^2$ i.e. $\sqrt{x^{2k-2h-1}} = u/v$.

Mais cela revient (quasiment) à démontrer l'irrationalité de la racine carrée de n'importe quelle puissance impaire de *x*,[116] autrement dit à prouver la question dont on est parti[117]. Nous sommes donc dans un cercle, et la généralisation échoue.

La démonstration proposée évite donc la critique aristotélicienne d'inadéquation du moyen à la fin que l'on peut soulever contre les démonstrations standard[118]. Elle n'encourt pas plus celles de bon sens, faites entre autres par G. Junge[119] et H. Vogt[120], qui tiennent pour invraisemblable qu'un intervalle de temps très long ait pu s'écouler entre la démonstration standard de l'irrationalité de √2, et la preuve générale[121].

---

[116] La situation pourrait sembler plus compliquée même que celle de départ, où il s'agissait (seulement) de montrer l'irrationalité de la seule racine carrée de *x*. En fait elles sont équivalentes d'après l'égalité 2 du a), à savoir : $a\sqrt{t} = \sqrt{a^2 t}$. Comme dans le cas *x* = 2, l'irrationalité de *x* équivaut à celle de toutes les puissances impaires de *x* (cf. *supra* IV.a).

[117] Cela n'est pas tout à fait exact car, cette fois, on sait que si elles sont rationnelles, elles doivent être représentées comme quotient, non pas de deux entiers quelconques, mais de deux entiers non multiples de *x*. C'est d'ailleurs cela qui permettait de conclure pour *x* = 2. Cela permettrait également de conclure, quoique pour les seuls entiers *x* premiers, à condition d'utiliser la proposition VII.30, selon laquelle si un entier premier divise un produit de deux nombres, alors il divise l'un d'entre eux. Mais on est alors dans une théorie de factorisation des nombres en facteurs premiers, tout à fait sans commune mesure avec les résultats que nous avons utilisés dans notre démonstration. En outre, cette théorie donne immédiatement le résultat dans toute sa généralité, comme nous l'avons vu en considérant la troisième démonstration d'irrationalité de √*2* (cf. *supra*, première partie, chap. IV et en particulier, la note 49).

[118] Cela ne signifie pas qu'avant l'élaboration complète de la théorie des irrationnels, attribuée par la tradition à Théétète, on n'ait pu résoudre certains cas particuliers ne relevant pas de notre démonstration. Toutefois, cela nécessite l'utilisation de résultats dont on n'a pas besoin pour la démonstration d'irrationalité que nous avons donnée aux chapitres IV et V.

On peut ainsi, relativement facilement, prouver géométriquement que si *x* et *n* sont des entiers quelconques et *r* le reste de la division par *x* de *n*, alors $n^2$ et $r^2$ sont simultanément divisibles (ou pas) par *x*. Algébriquement en écriture moderne, on peut le voir de la manière suivante :

soit $n = qx + r$ ; on a : $n^2 = (qx + r)^2 = q^2 x^2 + 2qxr + r^2$

et donc $r^2$ diffère de $n^2$ d'un multiple de *x*.

Les carrés de *n* et de *r* sont donc simultanément divisibles (ou pas) par *x*. Nous appellerons ce résultat la propriété (*).

On peut alors prouver, que, au moins pour les premiers termes de la suite des entiers étudiés par Théodore rapportée dans le *Théétète* de Platon, leurs racines carrées sont irrationnelles (cf. *supra*, note 49). En effet, d'après le résultat 1, il suffit, pour cela, de montrer l'impossibilité de l'égalité $x = u^2/v^2$, avec *u* et *v* entiers *non divisibles* par *x* (cf. *supra* note 115). Cette égalité s'écrit encore : $xv^2 = u^2$. Il suffit donc, pour que √*x* soit irrationnel, de montrer que si un entier est non divisible par *x*, son carré l'est aussi. D'après la propriété (*), il suffit de le prouver pour les seuls nombres inférieurs (strictement) à *x* (le reste de la division par *x* est toujours strictement plus petit que *x*).

Considérons ainsi le cas *x* = 3. Les seuls entiers inférieurs à *3* sont *1* et *2*, et leurs carrés ne sont pas des multiples de *3*. De même pour *x* = 5, les seuls entiers inférieurs à *5* sont *1, 2, 3* et *4* et leurs carrés ne sont pas des multiples de *5*. Ce processus se poursuit très facilement jusqu'à *11*, puis il devient de plus en plus fastidieux, mais pas impossible, au-delà.

On peut remarquer que cette méthode ne donne pas seulement les cas d'irrationalité des racines carrées, mais aussi bien ceux où celles-ci sont rationnelles. Ainsi soit *x* = 9, on doit considérer les carrés des entiers *1, 2, 3, ..., 8*. On constate que $3^2 = 9$ est non seulement un multiple de *9*, mais lui est égal, et donc la racine carrée de *9* est égale à *3*, elle est donc rationnelle.

Par contre, cette méthode ne permet pas d'obtenir le résultat général, puisqu'il s'agit de considérer les entiers l'un après l'autre, et qu'en outre, elle est de plus en plus compliquée à mesure que ceux-ci deviennent plus grands.

Certes, établir la propriété (*) suppose un développement de l'arithmétique grecque incompatible avec l'époque où nous pouvons situer notre démonstration d'irrationalité de racine carrée de *2*. Cela montre néanmoins que des résultats intermédiaires allant au-delà ont très probablement été obtenus entre cette période et celle où a été élaborée une théorie des irrationnels à la manière des livres VII et VIII des *Éléments*. Ainsi, il n'est en rien invraisemblable qu'une longue période de temps ait pu s'écouler entre la première preuve d'irrationalité de racine carrée de *2* et le développement de la théorie arithmétique qui permet de résoudre le problème dans toute sa généralité.

On peut également noter que cette méthode fournirait une éventuelle démonstration alternative pour la leçon d'irrationalité de Théodore dans le *Théétète* de Platon. Cela montre une fois encore, l'étroite relation existant entre les deux questions.

[119] *Cf.* « Wann haben die Griechen das Irrationale entdeckt ? », dans *Novae Symbolae Joachimicae*, 1907, p. 221-264.

[120] *Cf.* « Die Geometrie des Pythagoras », dans *Bibliotheca Mathematica*, 1908-1909, p. 15-54 et « Die Entstehungsgeschichte des Irrationalen nach Plato ... », *ibid.*, 1910, p. 97-155.

[121] Heath pense également que la première démonstration d'irrationalité de racine carrée de *2* est très ancienne, et que c'est bien celle visée par les textes des *Analytiques*. Mais parce qu'il admet qu'elle est donnée par la preuve standard, il se propose d'expliquer cette latence de la manière suivante. Les mathématiciens grecs auraient été pris à cette époque par

## VI. Un résultat vieux de 4 000 ans ?

Nous avons insisté au début de cet article, sur l'ancienneté des propriétés utilisées dans notre démonstration, qui remontent pour la plupart au moins à deux millénaires BCE.

Pourtant, il n'est aucune trace de ce résultat (en tant que théorème général de mathématique) avant celles de la Grèce classique.

Au vu de la rareté des textes de cette époque parvenus jusqu'à nous, cette absence n'infirme pas nécessairement qu'il ait pu être établi ailleurs. Néanmoins, le manque de toute référence dans les textes grecs que nous avons est certainement significatif. Les Grecs anciens, et Platon en particulier, si soucieux d'accorder aux Égyptiens anciens les découvertes les plus étonnantes (et pas seulement en mathématiques), ne leur réfèrent jamais cette question.

Le théorème (dit) de Pythagore, quant à lui, est indépendant de ce résultat, et n'intervient pas dans la démonstration. Celle-ci est une preuve de non-existence : il n'est pas de nombre (rationnel) de carré égal à *2*, pas plus que de bouc-cerf ou de solide régulier (convexe) à dix côtés[122].

Mais l'existence même, dans le cadre géométrique, d'une grandeur dont le carré est égal à *2* est une conséquence du théorème de Pythagore, ou plus exactement, d'un cas particulièrement simple de celui-ci. La démonstration en est donnée dans le *Ménon* de Platon, Socrate faisant apparaître sur un dessin, ou plutôt conduisant, comme par la main, un très jeune serviteur grec à la découverte du doublement d'un carré, par considération de sa diagonale (cf. *supra* IV.b). Cette existence est préliminaire à l'établissement de quelque résultat que ce soit le concernant. La séparation entre la définition et son existence est un point sur lequel insiste Aristote[123], aussi bien que l'impossibilité de dire quoi que ce soit d'une chose avant d'en avoir prouvé l'existence[124].

L'irrationalité de racine carrée de *2*, outre sa formulation comme incommensurabilité, suppose donc ce théorème[125]. Son attribution à Pythagore pose des difficultés, les témoignages étant très tardifs, postérieurs d'au moins un demi-millénaire, et peu sûrs[126]. Certains historiens pensent qu'il est beaucoup plus ancien, le faisant remonter aux Mésopotamiens, aux Égyptiens anciens ou aux Indiens[127]. Il est vrai que des cas numériques (ainsi le triangle de côtés 3, 4 et 5) se retrouvent dans de très anciens textes qui nous sont parvenus, la difficulté étant le niveau de généralité de ces résultats, mais aussi ce qu'il faut entendre par « démonstration mathématique »[128].

Le cas de l'irrationalité est moins disputé. Il est vrai qu'il est explicitement attribué à Pythagore ou aux pythagoriciens par Aristote ainsi que dans de nombreux textes grecs, jamais aux Égyptiens de l'antiquité. Et sa preuve, chez Platon, Aristote aussi bien que chez les auteurs postérieurs, est rapportée à la théorie du pair et de l'impair, qui est synonyme de la théorie des nombres pour les mathématiciens grecs[129]. C'est d'ailleurs sur cet accord des

---

l'étude d'autres questions de (plus) grande importance, telles la quadrature du cercle, la trisection de l'angle ou le doublement du cube (Heath 1956, I, p. 413). Quoi qu'il en soit, l'objection des auteurs allemands ne peut tenir que si l'on suppose qu'une partie des résultats du livre VII des *Éléments*, y compris la proposition 22, devait être connue pour obtenir l'irrationalité de racine carrée de *2*, ce qui n'est pas le cas dans la démonstration que nous avons proposée (cf. *supra*, chap. IV).

[122] Exemple favori de Leibniz, dont l'impossibilité n'apparaît pas immédiatement (Heath 1956, I, p. 145).

[123] *An. post.* I, 10 ; II, 7, 92b10-39.

[124] *Ibid.* II, 8, 93a16-21.

[125] Comme le souligne Heath, pour qui cela relève de l'évidence (Heath 1956, I, p. 351).

[126] Ainsi Diogène Laërce (VIII, 12) rapporte une anecdote de seconde main, sans se prononcer lui-même, et Proclus, dans son commentaire d'Euclide, paraît plutôt sceptique. Thomas Heath, dans une étude minutieuse, considère néanmoins qu'il n'y a pas de bonnes raisons pour douter, sur ce point, de la tradition (Heath 1956, I, p. 350-364).

[127] Pour une position critique concernant cette dernière hypothèse, *cf.* Heath 1956, I, p. 360-364 ; pour un point de vue opposé, *cf.* Kim Plofker, *Mathematics in India*, Princeton, 2008.

[128] *Cf.* p. ex. l'ouvrage d'Eleanor Robson, *Mathematics in Ancient Iraq: A Social History*, Princeton, 2008 (p. 28 *sq.*).

[129] Voir la longue analyse que fait Maurice Caveing sur cette question (Caveing 1998, III, p. 142-145).

historiens que s'appuie Heath pour défendre l'authenticité de l'attribution du théorème de Pythagore à celui-ci[130].

Mais l'argument essentiel nous paraît se trouver ailleurs. Certes, l'ancienneté des résultats sur le pair et l'impair, utilisés dans la preuve de l'irrationalité de la racine carrée de *2*, est bien connue des historiens de l'antiquité égyptienne et mésopotamienne. Il est même possible qu'égyptiens et mésopotamiens aient été en possession d'une forme plus ou moins primitive du théorème de Pythagore. En revanche, ce qui paraît complètement absent de leurs mathématiques, c'est l'utilisation de la méthode dite de démonstration par l'impossible, où il s'agit paradoxalement de prendre comme point de départ un énoncé faux, pour en conclure une proposition mathématique (vraie).
Dans la dernière partie, nous allons voir le lien étroit entretenu par la preuve d'irrationalité que nous proposons et la méthode générale de raisonnement par l'impossible, ainsi que la relation de celle-ci au discours philosophique.

---

[130] *Cf.* Heath 1956, I, p. 351.

# Troisième partie

## Sur l'origine du raisonnement par l'impossible

La démonstration d'irrationalité que nous avons donnée dans la deuxième partie rend convenablement compte des témoignages textuels portant sur l'origine des irrationnels, contrairement aux démonstrations standard. Mais elle présente un autre intérêt concernant l'origine du raisonnement par l'impossible. Elle explique, chez les Grecs anciens, le caractère paradigmatique de la démonstration géométrique d'incommensurabilité, en tant que modèle de démonstration par l'impossible. Cette propriété d'incommensurabilité était, certes, en elle-même, fondamentale. Mais, n'utilisant aucun résultat intermédiaire déjà obtenu par cette méthode, c'était surtout le premier théorème qu'aucune méthode directe ne permettait d'obtenir[131]. D'où la relation entre irrationnel et preuve par l'impossible, que l'on retrouve, en particulier, dans les ouvrages aristotéliciens. La valeur de cette forme démonstrative dérivait de l'importance accordée à la « découverte » de l'irrationalité[132].

On peut se demander si cette relation entre démonstration par l'impossible et preuve d'incommensurabilité est établie par Aristote, ou bien s'il ne fait que rapporter le point de vue des mathématiciens. Sa présentation en deux temps, dans la *Métaphysique,* où tous s'étonnent d'abord de l'incommensurabilité de la diagonale au côté du carré, alors que, dit ensuite Aristote, les mathématiciens, eux, seraient frappés de stupéfaction si elle devenait commensurable[133], permet de pencher pour le second cas.

Toutefois, la démonstration par l'impossible n'est pas cantonnée à la sphère mathématique. Elle permet plus généralement d'établir ou réfuter des thèses[134], et elle est massivement utilisée aussi bien dans les dialogues platoniciens que par Aristote[135], qui lui consacre une forme de syllogisme (*An. pr.* I, 44).

Cela conduit Arpad Szabó à formuler une double thèse. D'une part, la cause de la singularité des mathématiques grecques antiques, en tant que démonstratives et hypothético-déductives, est à rechercher dans cette forme de raisonnement[136]. D'autre part, celle-ci serait née chez les Éléates, le premier texte où se trouverait un tel raisonnement serait le *Poème* de Parménide. Et de conclure au fondement philosophique des mathématiques grecques antiques, donc des mathématiques telles que nous les entendons. Il reprend cette argumentation dans Szabó 1988, ajoutant qu'inversement Aristote n'a eu aucune influence sur les mathématiques de son époque[137].

---

[131] Pour une discussion sur ce point, cf. L. Löwenheim, « On making indirect proofs direct », *Scripta mathematica*, 12, 1946, p. 125-139.

[132] Lorsqu'un mathématicien moderne, tel G. Hardy, recherche des modèles pour faire comprendre la nature des « vrais » théorèmes mathématiques, c'est-à-dire à la fois « beaux » et « sérieux », le meilleur qu'il trouve est précisément l'irrationalité de la racine carrée de *2* (Hardy 1985, p. 31). Quant le mathématicien contemporain Vladimir Arnold veut montrer toute l'importance de la preuve par Newton de la non-intégrabilité algébrique de certaines courbes, « prototype des démonstrations ultérieures » d'Abel et de Liouville, c'est encore, comme le physicien anglais lui-même, à celle de l'irrationalité des racines carrées d'entiers qu'il se ramène (*Huygens et Barrow, Newton and Hooke*, Birkhaüser, 1990, p. 94).

[133] *Metaph.* A, 2, 983a15-19.

[134] Ainsi Le Blond 1938, p. v, remarque : « Les affirmation d'Aristote sont souvent ponctuées, dans la plupart de ses traités, par les termes ἀνάγκη, ἀναγκαῖος, ἀναγκαίως, qui en soulignent la rigueur et qui rappellent sans cesse que la science véritable comporte la connaissance du nécessaire. » Et un peu plus loin : « Il est naturel, puisque εὔλογος signifie souvent une égalité ou une identité mathématique, qu'il exprime aussi la cohérence d'un raisonnement, la *rigueur logique*. » (p. 40).

[135] On en trouve un exemple dans la *Métaphysique* avec la démonstration rejetant une intelligence divine qui serait seulement « en puissance ». Si tel était le cas, il y aurait quelque chose de supérieur à la perfection divine. Car alors il serait « logique de supposer que la continuité de la pensée est pour elle une charge pénible. Ensuite, il est clair qu'il y aura quelque autre chose plus noble que l'Intelligence, à savoir l'objet même de la pensée. (…) L'Intelligence suprême se pense donc elle-même, puisqu'elle est ce qu'il y a de plus excellent. » (*Metaph.* L, 9, 1074b30-34, trad. Tricot.) Un autre exemple concerne sa polémique contre les platoniciens. Ainsi, les critiquant, il réfute longuement qu'une unité puisse ne pas être additionnable avec une autre unité. Mais c'est pour conclure que jamais ceux-ci n'avaient soutenu une telle doctrine. Et d'ajouter : « Pourtant, l'inadditionnabilité absolue des unités résulte logiquement de leurs principes » (*ibid*., M, 7, 1081a37-39). Mais, comme « en vérité, elle n'est pas admissible » (*ibid*.), cela invalide leur thèse, à la manière d'un raisonnement par l'absurde en mathématique, où la conclusion permet de prouver la fausseté de l'hypothèse initialement posée.

Pour défendre son point de vue, A. Szabó s'appuie sur le discours de Parménide, tenant un raisonnement par l'impossible, dans l'ouvrage éponyme de Platon. Mais rien ne permet de douter que celui-ci mette dans la bouche de ses personnages des notions ou raisonnements, sans signifier nécessairement qu'ils en sont les auteurs. Si même, telle était l'intention de Platon, cela signifierait plutôt que Parménide aurait été le premier à l'introduire dans le discours philosophique[138].

Pourtant, Wilbur Knorr remarquait déjà que, lorsqu'une telle direction existe, son sens en est généralement opposé, les philosophes s'emparant des travaux, résultats, définitions ou démonstrations mathématiques, afin de les exploiter[139]. Et dans l'*Euthydème*, ce sont les dialecticiens que Socrate voit « cuisiner » les résultats mathématiques[140]. C'est pourquoi les philosophes doivent être rompus aux mathématiques (cf. *Resp.* VII)[141], ce que veut sans doute également signifier la tradition autour de l'inscription sur le portique de l'Académie[142].

Que ce soit chez Platon ou chez Aristote, voire chez les commentateurs plus tardifs, on ne trouve pas d'allusion à une origine philosophique, moins encore parménidienne, du raisonnement par l'impossible. Si cette relation existait, ce silence serait difficilement explicable, au vu de l'importance accordée, chez ces auteurs, au rapport entre mathématique et philosophie.

En outre, la tradition rapporte que nombre de penseurs présocratiques ont été mathématiciens, ainsi Thalès, Démocrite, Anaximandre ou Pythagore. On peut donc douter qu'une nouvelle forme de raisonnement eût pu être utilisée longtemps dans l'un des champs, sans passer dans l'autre. Celui par l'impossible fut très probablement mis en œuvre, antérieurement à Parménide, concomitamment ou de manière très proche, en mathématique et dans les autres secteurs de connaissances à la recherche de la rigueur, les preuves mathématiques, particulièrement celles obtenues par l'impossible, étant de l'ordre de l'évidence[143].

Il n'empêche que l'on peut suivre Arpad Szabó, lorsqu'il propose de lier naissance des mathématiques grecques et raisonnement par l'impossible. Et il est vrai que beaucoup d'historiens des sciences considèrent que les mathématiques qui précèdent celles de la Grèce antique, quelque développées qu'elles aient été, formaient des techniques caractérisées par un aspect purement utilitaire[144]. Sans aborder ici la question de l'exceptionnalité de la mathématique grecque[145], Szabó a certainement raison de souligner l'importance de ce type de raisonnement en mathématiques, pour dépasser la représentation sensible, en particulier visuelle[146]. Cette méthode serait à l'origine de la construction démonstrative, au sens hypothético-déductif, telle qu'elle apparaît, pour autant qu'on le sache, en Grèce antique. L'introduction des grandeurs irrationnelles, en mettant fin aux seuls calculs sur les entiers, exigeait cette nouvelle approche, ainsi dans la preuve de l'irrationalité de racine carrée de *2*.

---

[136] *Cf.* Szabó 1977, troisième partie, notam. p. 336-337.

[137] Szabó 1988, p. 287, 290-291.

[138] Et Zénon, en aurait fait de même pour la contraposition (cf. *supra*, note 71).

[139] *Cf.* Knorr 1982, p. 121, 135, 143. Granger 1998, p. 23, trouve également douteuses les conclusions de Szabó. Même opinion chez J. Vuillemin concernant le sens des rapports généraux entre mathématiques et philosophie. Pour lui, les nouveautés mathématiques se répercutent sur la philosophie (*La Philosophie de l'algèbre*, p. 4). Voir aussi Bourbaki 1960, p. 11.

[140] *Euthydème*, 290c.

[141] Selon Gregory Vlastos, le niveau de connaissance de Platon serait l'équivalent, à son époque, de celui exigé de nos jours pour un doctorat (Vlastos 1992, p. 138). Il s'agit pour les philosophes d'utiliser les mathématiques, non d'un appel à abandonner la philosophie à leur profit, comme l'a fait Théodore, le géomètre du *Théétète* (165a). Reproche qu'Aristote adresse aux philosophes de son temps (*Métaphysique*, I, 9, 992a35-b1).

[142] 'Ἀγεωμέτρητος μηδεὶς εἰσίτω' (« Que nul n'entre ici qui ne soit géomètre »). Sur son caractère certainement apocryphe, *cf.* Saffrey 1968 et aussi Ofman, à paraître, Annexe 2.

[143] *Cf.* Arist. *An. pr.* I, 44, 50a39 et *supra*, première partie, chap. V.

[144] Ainsi L. Robin écrit : « Jamais, que nous sachions, la science orientale, à travers tant de siècles d'existence, et même après qu'elle eut pris contact avec la science des Grecs, ne paraît avoir dépassé les préoccupations utilitaires ou les curiosités de détail, pour s'élever à la pure spéculation et à la détermination des principes. » (Robin 1973, p. 51-52). Voir aussi Gaston

Cela est remarqué par de nombreux auteurs. Ainsi B. Artmann remarque que si, dans la justification d'une démonstration où se trouvent des constructions géométriques, celles-ci peuvent ou non aider,

> *la situation change complètement avec l'introduction de concepts abstraits qui n'ont pas de contrepartie intuitive. Historiquement, le premier et le plus important d'entre eux est le concept d'incommensurabilité, ou d'irrationalité. Il n'y a pas d'arrière-plan intuitif qui permettrait de contrôler la validité des assertions concernant de tels concepts abstraits*[147].

W. Knorr note que « contrairement aux autres branches des mathématiques, ce domaine ne peut faire appel à la pratique pour soutenir ses affirmations. (…) Sa justification ne peut être basée que sur une *déduction logique*, et une théorie de telles quantités ne peut avoir qu'*un* critère de validité : *la cohérence logique.* »[148] (nous soulignons).

Et Maurice Caveing considère, plus généralement, qu'avec le raisonnement par l'absurde, la géométrie est en quelque sorte contrainte de sortir du cadre de l'image (Caveing 1998, III, p. 316). C'est pourquoi l'irrationalité se retrouve au centre de la démonstration par l'impossible[149]. Ce serait alors naturellement qu'elle passerait des mathématiques au discours fondé sur le raisonnement démonstratif.

Les mathématiques fournissaient en effet le premier type de discours cohérent, aussi bien dans sa totalité que dans ses parties, ce qui, selon Alexander Nehamas, représente le caractère opératoire distinguant, chez Platon, la connaissance (ἐπιστήμη) de la croyance fausse[150]. Et G.G. Granger caractérise la « connaissance rationnelle [comme celle qui] conçoit et organise son objet de telle sorte que ses parties ne soient pas reconnues comme *contradictoires* entre elles, ni contradictoires avec le tout »[151]. Les mathématiques apparaissant alors comme un discours « absolument » rationnel[152].

---

Milhaud, *Nouvelles études sur l'histoire de la pensée scientifique*, 1911, p. 45 ; Abel Rey : *La Science orientale avant les Grecs*, 1929, p. 271 ; P.-H. Michel, *De Pythagore à Euclide*, 1950, p. 12. Pour un point de vue différent, voir Karpinsky 1937 ; Bourbaki 1960, p. 9-10. Voir aussi *supra*, note 126.

[145] Pour certains, la science elle-même commence avec les Grecs de l'époque classique. Ainsi John Burnet considère que la « science naturelle est une création des Grecs » et qu'on ne trouve « pas la moindre trace de science en Égypte ni même à Babylone » (*Greek Philosophy*, I, p. 4-5).

[146] Cf. *Resp.* VII, 534c. M. Caveing va plus loin, considérant que l'espace euclidien s'oppose à l'espace visuel comme le monde sphérique au monde elliptique (Caveing 1988, p. 294).

[147] 'What has been said so far is about geometric construction, where pictures may or may not help. The situation changes completely with the introduction of abstract concepts that have no intuitive counterpart. Historically, the first and most significant one of these is the concept of incommensurability, or irrationality. There is no intuitive background that would allow one to control the validity of assertions about such abstract concepts.' (Artmann 1999, p. 111.)

[148] 'Unlike other branches of mathematics, this field can make no appeal to practice for supporting its claims. Many mathematical traditions confronted the computational anomaly posed by such quantities as √2; there was little problem in producing approximations which would satisfy any given context. But the Geeks advanced to the assertion that such quantities were irrational. This is a statement of an entirely different theoretical order. Its justification can be based only upon logical deduction *and a theory of such quantities can have* but one *criterion of validity:* logical consistency.' (Knorr 1975, p. 4.)

[149] Au tout début du livre A de la *Métaphysique* (2, 983a13-20), « l'incommensurabilité de la diagonale avec le côté du carré » apparaît comme le modèle de la connaissance scientifique. Dans les *Premiers analytiques*, c'est comme exemple d'un raisonnement par l'absurde que, par trois fois, Aristote se réfère à la preuve d'incommensurabilité. Ce sont précisément ces textes qui nous permettent, aujourd'hui, d'en tenter une reconstruction, ainsi que nous l'avons fait dans la deuxième partie.

[150] « Nous possédons l'*episteme*, quand nous maîtrisons la structure axiomatique du système en question et quand nous pouvons prouver la vérité de n'importe lequel de ses éléments. » (Nehamas 1991, p. 286.)

[151] Granger 1976, p. 270-271 ; c'est l'auteur qui souligne.

[152] Au sens de « l'absolue continuité » dont Socrate qualifie la « théorie » d'Hippias selon laquelle pour qu'une chose ait une propriété il faut (et il suffit) que chacune de ses parties la possède (*Hipp. maj.* 300b-302c). Socrate la réfute de cette

Cette démarche, caractérisée par la recherche de la certitude vraie opposée à la vraisemblance, est considérée comme un trait originel des mathématiques. Ainsi, dans son commentaire du livre I des *Éléments* d'Euclide, plus particulièrement du postulat I.5 (dit des parallèles), Proclus approuve un auteur plus ancien, Géminus, d'affirmer que

> *nous avons appris des pionniers mêmes de cette science [la géométrie] de n'avoir pas le moindre égard pour les images vraisemblables lorsqu'il s'agit d'une question concernant les raisonnements à inclure dans notre doctrine géométrique*[153].

De ce point de vue, les mathématiques pouvaient apparaître comme un modèle à imiter.

Toutefois, qu'un tel discours soit possible, hors de la sphère mathématique, n'est pas une évidence. Les énoncés usuels sont des conjectures, des possibilités ou des probabilités bien plutôt que des certitudes[154].

Aristote lui-même, lorsqu'il polémique contre les Mégariques, est confronté à des problèmes soulevés dans le cadre d'une logique inusuelle, admettant diverses valeurs de vérité (le nécessaire et l'impossible, mais aussi le possible ou le contingent, le probable …). Il en est de même à propos des énoncés concernant le futur, l'exemple le plus fameux étant celui sur lequel ont travaillé des générations entières de logiciens médiévaux, d'une bataille navale pour le lendemain[155]. Selon le Stagirite, ces énoncés participent d'une logique de ce qui est en puissance opposée à celle de ce qui est en acte, où le principe du tiers exclu (une proposition doit être ou bien vraie, ou bien fausse) peut être tenu en échec[156]. Car alors, ni l'affirmation que demain une bataille navale aura lieu ni sa contradictoire, « une bataille navale n'aura pas lieu », ne sont nécessaires. Par contre, elles ne peuvent être simultanément vraies[157].

Le principe commun au raisonnement admettant la réfutation, n'est donc pas, comme en mathématiques, telles que les concevaient du moins les mathématiciens de la Grèce antique, celui du tiers exclu, mais un principe plus faible, celui de non-contradiction.

Suivant celui-ci, un discours contenant deux énoncés contradictoires est nécessairement faux. Or, si l'on ne cherche qu'à *réfuter* une thèse (et non pas à en *établir* une), la démonstration par l'impossible ne nécessite que cette forme faible. On comprend alors l'importance qu'Aristote lui accorde, ce principe apparaissant bien plus fréquemment dans ses textes que celui du tiers exclu. Et là même où il s'applique, à l'inverse d'un raisonnement mathématique où l'énoncé faux, à savoir l'hypothèse de départ, est évident, la source de l'erreur dans un discours général est souvent plus difficile à discerner[158]. Le raisonnement

---

manière : « il n'est donc pas absolument nécessaire, comme tu le disais à l'instant, que ce que nous sommes tous les deux ensemble, chacun des deux le soit aussi, ni non plus que ce qu'est chacun, tous les deux le soient ensemble. » (*Hipp. maj.* 302b1-2.)

[153] '… *that we have learned from the very pioneers of this science not to have any regard to mere plausible imaginings when it is a question of the reasonings to be included in our geometrical doctrine.*' Et d'ajouter l'opposition entre rhétorique et géométrie soulignée par Aristote, notant qu'il serait aussi absurde « de demander des démonstrations scientifiques à un rhétoricien que d'accepter de simples vraisemblances d'un géomètre ». ('*For Aristotle says that it is as justifiable to ask scientific proofs of a rhetorician as to accept mere plausibilities from a geometer*') (Morrow 1992, p. 202).

[154] *Cf.* p. ex. Knorr 1975, p. 14 ; Pellegrin 2007, p. 35.

[155] Cf. *Metaph.* Q ; *De int.* 9 ; voir aussi Granger 1998, p. 224-225 et Gourinat 2001. Pour un commentaire littéral de ce texte, cf. par exemple Anscombe 1967.

[156] *De int.* 9, 18b32-39.

[157] *Ibid.*, 18b17-25, 19b30-31. Ce que Quine critique en qualifiant cette position d'assez désespérée pour oser penser qu'une disjonction (*p* ou *q*) puisse être vraie sans que l'un des termes (soit *p*, soit *q*) le soit (« On a supposed Antinomy », dans *The Ways of Paradox and Other Essays*, Random House, 1966, p. 21 ; cf. aussi l'article de J. Vuillemin in Vuillemin 2001, p. 78-87). Sur la question du déterminisme ou de l'indéterminisme aristotélicien pour expliquer ce passage, cf. par exemple Sorabji 1983, chap. 5 (« Tomorrow's Sea Battle: an argument from past truth ». J. Hintikka a proposé une interprétation très contestée sauvegardant le tiers exclu (Hintikka 1964). Mais c'est au prix d'une logique temporelle qui est différente de la logique usuelle (*ibid.*, p. 475). Renforçant encore la différence entre philosophie et mathématiques, elle ne contrevient en rien à l'argumentation que nous développons ici.

[158] D'où cette difficulté d'interpréter les dialogues de Platon, certains allant jusqu'à se demander si le philosophe athénien n'a pas pris « un malin plaisir à crypter » ses œuvres (*cf.* M. Crubellier-P. Pellegrin, « Approches de la Physique d'Aristote », *Oriens-Occidens*, 2, 1998, p. 1-37, note 8 concernant le *Timée*). Lorsqu'on aboutit à une aporie, une hypothèse de départ est nécessairement fausse, mais laquelle ? Ainsi, à propos du « jugement faux qui fait croire qu'un jugement faux est

discursif général n'est donc pas restreint aux seuls cas où s'applique le principe du tiers exclu, et par conséquent, au domaine d'application du raisonnement par l'absurde, qui, comme on l'a vu, permet d'établir des propositions vraies, et non pas (seulement) la fausseté de tel ou tel énoncé[159]. Il est ainsi extrêmement improbable que celui-ci soit apparu ailleurs que là où il s'imposait naturellement, et le principe du tiers exclu assuré, au sein des mathématique, pour lesquelles ce type de raisonnement est parfaitement adapté, et simple à mettre en œuvre [160]. Les démonstrations s'y fondent en effet sur une logique strictement bivalente du tiers exclu (un énoncé est ou vrai ou faux, « ou » étant exclusif).

À l'inverse, le raisonnement par l'impossible n'a rien de naturel à l'intérieur du discours philosophique, ce que montrent par exemple les énoncés aristotéliciens concernant les futurs contingents[161]. Il a donc été très certainement importé par ceux qui tenaient à donner un caractère rigoureux au discours. Cela a été d'autant plus facile et rapide que, à l'origine, les divers champs de connaissances n'étaient pas disjoints[162]. Par la suite, si les savants n'étaient pas savants en tout (comme se présente encore Hippias dans les dialogues éponymes de Platon), du moins connaissaient-ils bien ce qui se faisait dans les autres domaines. Ainsi Platon nous montre Socrate, mais aussi Parménide, excellents connaisseurs des mathématiques[163]. Et Abel Rey ne devait pas être si éloigné de la vérité lorsque, sur le lien entre démonstration par l'impossible en mathématique et en philosophie, il écrivait :

> *Est-ce trop s'aventurer que supposer que ce schéma [l'argumentation par l'impossible] pouvait être d'usage courant et traditionnel, dans l'argumentation, au moment où l'emploie Zénon ? La démonstration par l'absurde sera affectionnée de la mathématique grecque. Zénon n'a fait que la transposer à la logique pure* (Rey 1933, p. 202).

Si l'on devait choisir, il est donc bien plus vraisemblable que, contre l'opinion de Szabó, l'origine de la démonstration par l'impossible se trouve en mathématiques[164].

---

impossible » dans le *Théétète*, Burnyeat 1990, p. 94, analyse de manière positive la suggestion d'E. Zeller, qui propose de voir dans la page 201 « une sorte de *reductio ad absurdum* ». Il ajoute néanmoins qu'« il vaut mieux être circonspect à propos de cette suggestion, car Platon *nous a laissé le soin de localiser* la supposition erronée » (nous soulignons). Pour des exemples de tels raisonnements en philosophie, et la complexité à comprendre leurs conclusions, *cf. supra*, note 133.

[159] Cela fait sens avec l'analyse que fait J.-M. Le Blond de l'emploi par Aristote du terme εὔλογος sous une forme négative. Il l'utilise toujours pour désigner l'incohérence (ce qui n'est pas εὔλογος) que l'on obtient par cette déduction. Comme on l'a vu (cf. *supra* note 132), ce terme est suffisamment vague pour recouvrir aussi bien le sens d'erreur logique que d'invraisemblance, au sens de ce qui est improbable. Le Stagirite semble vouloir distinguer ainsi, dans les réfutations philosophiques, l'impossibilité absolue de type mathématique (ce qui va à l'encontre du nécessaire), de ce qui est improbable (ce qui va à l'encontre de l'opinion de tous ou du plus grand nombre).

[160] *Cf.* Knorr 1975, p. 3.

[161] Mais aussi la possibilité de dissocier cohérence locale et globale. Un exemple en est donné par Moravcsik 1967, p. 10, lorsqu'il défend une approche pluraliste, par morceaux, de la philosophie d'Aristote ('*the proper approach to Aristotle's philosophy is a piecemeal or pluralistic one.*'). En effet, écrit-il, « ce qu'Aristote dit comme réponse à une question peut ne pas être lié à ses réponses à d'autres questions en apparence similaires. Aussi, *accepter certaines suggestions d'Aristote ne conduit pas à accepter la totalité de ses affirmations philosophiques*. Cette remarque s'applique aussi bien au contenu qu'à la méthodologie. » ('*What Aristotle says as an answer to one question may not relate to his answers to other questions similar in appearance. Thus acceptance of some of Aristotle's suggestions does not commit one to a wholesale acceptance of his philosophic claims. This remark applies both to content and to methodology.*' (Nous soulignons.)

[162] Selon Heath 1981, p. 3, le génie des Grecs anciens en mathématiques « était simplement un aspect de leur génie pour la philosophie. Leurs mathématiques constituaient en fait une large partie de leur philosophie jusqu'à Platon. Toutes deux avaient la même origine." ('*Their genius [of the Greeks] for mathematics was simply one aspect of their genius for philosophy. Their mathematics indeed constituted a large part of their philosophy down to Plato. Both had the same origin.*')

[163] Knorr 1986 note ce qu'il appelle la « convergence des intérêts géométriques et philosophiques » à l'époque de Platon (p. 87, 111). Selon lui, une conséquence malheureuse en histoire des sciences a été d'amener les spécialistes à déformer « l'interprétation de la géométrie ancienne » en fonction de leur intérêt (p. 87).

[164] Celui-ci examine cette possibilité (Szabó 1977, p. 269-273), et la récuse d'après des arguments de datation, dont il reconnaît, lui-même, le caractère spéculatif (ibid. p. 271). Dans le cadre de l'interprétation des textes aristotéliciens, aussi bien Barnes 1981, p. 18-19 que Pellegrin 2005, p. 48-49, soulignent la difficulté de décider de la direction de l'influence entre Aristote et les mathématiques.

Pour conclure, il est difficile de ne pas ajouter une pièce aux conjectures déjà existantes concernant la période où a été obtenue l'irrationalité de racine carrée de *2* (ou plus précisément, l'incommensurabilité du rapport de la diagonale aux côtés du carré).

Suivant notre démonstration, les résultats mathématiques nécessaires étaient très anciennement disponibles[165]. Reste la mise en forme du raisonnement qui exige une démonstration par l'impossible (cf. *supra*, début de ce paragraphe). Ou en d'autres termes, une conception nouvelle de la mathématique, permettant d'élaborer cette forme de preuve. Cela irait dans le sens de, ou pour être plus exact, ôterait certains obstacles à, une datation haute. Elle pourrait en effet remonter à la période pythagoricienne, soit, si l'on suit la tradition, le début du v$^e$ siècle, voire antérieurement. Une difficulté majeure est de dater cette forme de raisonnement, la circularité guettant toujours : d'une datation hypothétique de celle-ci, on tirerait celle-là, et inversement[166]. Quoi qu'il en soit, suivant la perspective que nous proposons, l'irrationalité se trouve intrinsèquement liée à la méthode par l'impossible, non seulement comme un cas d'application fondamental, mais logiquement, voire chronologiquement, en accord avec la stricte relation que l'on trouve dans les textes.
Ce n'est pas le lieu de procéder ici à une critique détaillée des arguments soulevés par Knorr, contre une datation haute. Cette position s'appuie sur la thèse selon laquelle une théorie générale de l'irrationalité aurait été développée par Théétète. Pourtant, argumentant d'après les mêmes lignes, Thomas Heath aboutit à une telle datation (*cf.* Heath 1981, p. 156-157). Suivant notre interprétation, les deux questions sont indépendantes, au moins du point de vue de leur démonstration. Car il est bien une relation entre elles : la nécessaire antériorité de l'incommensurabilité de la diagonale et de ses conséquences (*cf.* chap. IV de la deuxième partie) par rapport à la preuve générale[167].

---

[165] Le plus récent serait sans doute une certaine forme du théorème de Pythagore, donnant la diagonale du carré en fonction de ses côtés. Non qu'il participe de la preuve, comme on peut s'en convaincre par celles que nous avons données dans les première et deuxième parties, mais pour une question d'existence. En effet, sans lui, il n'y aurait pas de raison de seulement penser l'incommensurabilité (*cf.* Platon, *Lois*, VII, 820a-c).

[166] *Cf.* aussi *supra* note 29.

[167] Ceci sera développé et approfondi dans notre article sur le *Théétète* à paraître dans un prochain numéro.

# Ouvrages cités

Allen 2000 : James Allen, *Middle Egyptian. An introduction to the Language and Culture of Hieroglyphs*, Cambridge Univ. Press, 2000

Anscombe 1967 : Elizabeth Anscombe, Aristotle and the Sea Battle, *De Interpretatione*, Chapter IX, in *Aristotle, a Collection of Critical Essays*, dir. J. Moravcsik, Doubleday Anchor, 1967, p. 15-33

Tricot 1991 : Aristote, *Métaphysique*, trad. J. Tricot, Vrin, 2 vol., 1991

Tricot 1992 : Aristote, *Premiers analytiques*, trad. J. Tricot, Vrin, 1992

Tredennick 1938 : Aristotle, *Prior Analytics*, trad. H. Tredennick, Harvard Univ. press, 1983 (1938)

Pellegrin 2005 : Aristote, *Seconds Analytiques*, trad. P. Pellegrin, Flammarion, 2005

Dalimier-Pellegrin 2004 : Aristote, *Traité du ciel*, trad. C. Dalimier-P. Pellegrin Flammarion, 2004

Pellegrin 1993 : Aristote, *Les Politiques*, trad. P. Pellegrin, Flammarion 1993

Pellegrin 2007 : Aristote, *Catégories – Sur l'interprétation (Organon I-II)*, trad. M. Crubellier, C. Dalimier, P. Pellegrin, Flammarion 2007

Aubenque 1962 : Pierre Aubenque, *Le problème de l'être chez Aristote*, PUF, 1962

Artmann 1999 : Benno Artmann, *Euclid, the Creation of Mathematics*, Springer, 1999

Barnes 1981 : Jonathan Barnes, Proof and the syllogism, in *Aristotle on Science:* The Posterior Analytics, E. Berti (dir.), Antenore, 1981, p. 17-59

Becker-Hofmann : Oskar Becker-Joseph Hofmann, *Histoire des mathématiques*, trad. R. Jouan, Lamarre, 1956

Bourbaki 1960 : Nicolas Bourbaki, *Éléments d'histoire des mathématiques, Hermann*, 1960

Burnyeat 1990 : Myles Burnyeat, Introduction au *Théétète* de Platon, trad. M. Narcy, Presses Universitaires de France, 1990

Caveing 1998 : *La constitution du type mathématique de l'idéalité dans la pensée grecque*, 3 vol., Presse Univ. du Septentrion, 1997

Caveing 1982 : *Zénon d'Élée, prolégomènes aux doctrines du continu*, Vrin, 1982

Caveing 1988 : Aristote et les mathématiques de son temps, in *Aristote aujourd'hui*, dir. M. Sinaceur, Érès, 1988, p. 293-299

Vitrac 19••-19•• : Euclide, *Éléments*, trad. B. Vitrac, 4 vol., PUF, 1994- ?

Goldschmidt 2003 : Victor Goldschmidt, *Le paradigme dans la dialectique platonicienne*, Vrin, 2003

Gourinat 2001 : Jean-Baptiste Gourinat, Principe de contradiction, principe du tiers-exclu et principe de bivalence : philosohie première ou organon ?, *Logique et métaphysique dans l'Organon d'Aristote*, M. Bastit-J. Follon ed., Peeters, 2001, p. 63-91

Granger 1998 : Gilles-Gaston Granger, *L'irrationnel*, Odile Jacob, 1998

Granger 1976 : *La théorie aristotélicienne de la science,* Aubier Montaigne, 1976

Hardy 1985 : Godfrey Hardy, *L'apologie d'un mathématicien*, trad. D. Julien-S. Yoccoz, Belin, 1985

Heath 1956 : Thomas Heath, *The thirteen Books of Euclid's Elements*, 3 vol., Dover, 1956 (1925)

Heath 1998 : *Mathematics in Aristotle*, Thoemmes Press, 1998 (1948)

Heath 1981 : *A History of Greek Mathematics*, *From Thales to Euclid*, Dover, 1981 (1921)

Hintikka 1964 : Jaakko Hintikka, The once and Future Sea Fight: Aristotle's Discussion of Future Contingents in *De Interpratione* IX, *Philosophical Review*, 1964, p. 461-492

Ifrah 1994 : Georges Ifrah, *Histoire universelle des chiffres*, Laffont, 1994 (1981), 2 vol.

Karpinsky 1937 : Louis Charles Karpinsky, Is there progress in mathematical discovery and did the Greeks have analytic geometry?, *Isis*, 27, 1937, p. 46-52

Kline 1982 : Morris Kline, *Mathematics, The Loss of Certainty*, Oxford University Press, 1982

Knorr 1982 : Wilbur Knorr, Infinity and Continuity: The Interaction of Mathematics and Philosophy in Antiquity, in *Infinity and Continuity*, dir. Kretzmann, p. 112-145

Knorr 1975 : *The Evolution of the Euclidean Element*, Reidel, 1975

Knorr 1986 : *The Ancient Tradition of Geometric Problems*, Birkhäuser, 1986

Le Blond 1938 : Jean-Marie Le Blond, ΕΥΛΟΓΟΣ *et l'argument de convenance chez Aristote*, Belles lettres, 1938

Moravcsik 1967 : Julius Moravcsik, Introduction, in *Aristotle, a Collection of Critical Essays*, Doubleday Anchor, 1967, p. 1-11

Nehamas 1991 : Alexander Nehamas, Le paradoxe de Ménon et Socrate dans le rôle d'enseignant, in *Les paradoxes de la connaissance, essais sur le* Ménon *de Platon*, dir. M. Canto-Sperber, Odile Jacob, 1991, p. 271-298, repris de *Oxford Studies in Ancient Philosophy*, 3, 1985, 1-30

Bertier 1978 : Nicomaque de Gérase, *Introduction arithmétique*, trad. J. Bertier, Vin, 1978

Ofman (à paraître) : Salomon Ofman, *Mouvement et origine de l'infinitésimal : Aristote, Euclide, Galilée*, à paraître

Penrose 2007 : Roger Penrose, À la découverte des lois de l'univers, trad. C. Laroché, Odile Jacob, 2007 (2004)

Dixsaut 1991 : Platon, *Phédon*, trad. M. Dixsaut, Flammarion, 1991

Brisson 1999 : Platon, *Parménide*, trad. L. Brisson, Flammarion, 1999

Pradeau 2005 : Platon, *Hippias Majeur, Hippias Mineur*, J.F. Pradeau, Flammarion, 2005

Brisson & Pradeau : Platon, *Les Lois*, trad. L. Brisson et J.F. Pradeau, Flammarion, 2 vol., 2006

Canto 1987 : Platon, *Gorgias*, trad. M. Canto, Flammarion, 1987

Croiset 1923 : Platon, *Gorgias, Ménon*, trad. A. Croiset, Belles Lettres, 1974 (1923)

Lamb 1925 : Platon, *Lysis, Symposium, Gorgias*, trad. W. Lamb, Harvard Univ. Press, 2001 (1925)

Morrow 1992 : Proclus, *A commentary on the first book of Euclid's* Elements, trad. G. Morrow, Princeton University Press, 1992

Rey 1933 : Abel Rey, *La jeunesse de la science grecque*, Renaissance du livre, 1933

Robin 1973 : Léon Robin, *La pensée grecque et les origines de l'esprit scientifique*, Albin Michel, 1973 (1923)

Russell 1970 : Bertrand Russel, Introduction à la philosophie mathématique, trad. G. Moreau, Payot, 1970

Saffrey 1968 : Henri-Dominique Saffrey, ΑΓΕΩΜΕΤΡΗΤΟΣ ΜΗΔΕΙΣ ΕΙΣΙΤΩ. Une inscription légendaire, *Revue des études grecques*, 81, 1968, p. 67-87

Sorabji 1983 : Richard Sorabji, *Necessity, Cause and Blame. Perspectives on Aristotle's theory*, Cornell Univ. Press, 1983

Szabó 1977 : Árpád Szabó, *Les débuts des mathématiques grecques*, tr. M. Federspiel, Vrin, 1977

Szabó 1988 : Árpád Szabó, Quelques problèmes fondamentaux des mathématiques grecques au temps d'Aristote et la prise de position du Stagirite, in *Aristote aujourd'hui*, dir. M. Sinaceur, Erès, 1988, p. 276-292

Vlastos 1992 : Gregory Vlastos, Elenchus and Mathematics: A Turning-Point in Plato's Philosophical Development, in *Essays on the Philosophy of Socrates*, dir. H. Benson, Oxford Univ. Press, 1992, p. 137-161

Vuillemin 2001 : Jules Vuillemin, *Mathématiques pythagoriciennes et platoniciennes*, Blanchard, 2001

Wedberg 1955 : Anders Wedberg, *Plato's Philosophy of Mathematics*, Almquist et Wiksell, 1955